\documentclass[12pt]{amsart}

\usepackage{pstricks}
\usepackage{pstcol}
\usepackage{graphicx}
\usepackage{amsmath}
\usepackage{amsfonts}
\usepackage{amssymb}
\usepackage{latexsym}

 \advance\hoffset by -0.4 truecm

\newtheorem{theorem}{Theorem}
\newtheorem{corollary}{Corollary}

\newcommand{\s}{\vspace{0.2cm}}
\newcommand{\K}{\text{\bf K}}
\newcommand{\R}{\text{\rm Re }}

\newcommand{\Diff}{\text{\rm Diff }}
\newcommand{\Rot}{\text{\rm Rot }}
\newcommand{\Vect}{\text{\rm Vect }}

\newcommand{\ess}{\text{\rm ess }}
\def\Vb{\text{\bf V}}

\begin{document}
\title[Evolution dynamics of conformal maps with quasiconformal extensions]
{Evolution dynamics of conformal maps with quasiconformal
extensions}
\thanks{This work is partially  supported by
Projects Fondecyt (Chile) \# 1030373, \#1040333, and UTFSM
\#12.03.23}
\author[Alexander Vasil'ev]{Alexander Vasil'ev}
\address{Departamento de Matem\'atica, Universidad T\'ecnica
Federico Santa Mar\'{\i}a, Casilla \newline 110-V,
Valpara\'{\i}so, Chile} \email{alexander.vasiliev@mat.utfsm.cl}
\begin{abstract}
We study one-parameter curves on the universal Teichm\"uller space
$T$ and on the homogeneous space $M=\Diff S^1/\Rot S^1$ embedded
into $T$. As a result, we deduce evolution equations for conformal
maps that admit quasiconformal extensions and, in particular, such
that the associated quasidisks are bounded by smooth Jordan
curves. Some applications to Hele-Shaw flows of viscous fluids are
given.
\end{abstract}

\maketitle

\section{Introduction}

Let $U$ denote the unit disk in the Riemann sphere $\hat{\mathbb C
}$ and $U^*=\hat{\mathbb C }\setminus \hat{U}$, where $\hat U$ is
the closure of $U$. By $S$ we denote the class of all holomorphic
univalent functions in $U$ normalized by
$f(\zeta)=\zeta+a_2\zeta^2+\dots$, $\zeta \in U$, and by $\Sigma$,
the class of all univalent meromorphic functions in $U^*$
normalized by
$\displaystyle{f(\zeta)=\zeta+c_0+\frac{c_1}{\zeta}+\dots}$,
$\zeta \in U^*$, $\Sigma_0$ stands for all functions from $\Sigma$
with $c_0=0$. These  classes have been one of the principal
subjects of research in Complex Analysis for a long time. The most
inquisitive  problem for the class $S$ posed by L.~Bieberbach in
1916 \cite{Bib} finally has been solved in 1984 by L.~de~Branges
\cite{Branges} who proved that $|a_n|\leq n$ for any $f\in S$ and
the equality is attained only for the Koebe function
$k(z)=z(1-ze^{i\theta})^{-2}$, $\theta\in [0,2\pi)$. The main tool
of the proof turned out to be a parametric representation of a
function from $S$ by the L\"owner homotopic deformation of the
identity map given by the {\it L\"owner differential equation}.
The parametric method  emerged almost 80 years ago in the
celebrated paper by K.~L\"owner \cite{Loewner} who studied a
one-parameter semigroup of conformal one-slit maps of $U$ coming
then at an evolution equation called after him. His main
achievement was an infinitesimal description of a semi-flow of
such maps by the Schwarz kernel that led him to the L\"owner
equation. This crucial result was generalized, then, in several
ways. Attempts have been made to derive an equation that allowed
to describe  a representation of the whole class $S$. Nowadays, it
is rather difficult to follow the correct history line of the
development of the parametric method because in the middle of the
20-th century a number of works dedicated to this general equation
appeared independently. In particular, P.~P.~Kufarev
\cite{Kufarev} studied a one-parameter family of domains
$\Omega(t)$, and regular functions $f(z,t)$ defined in
$\Omega(t)$. He proved differentiability of $f(z,t)$ with respect
to $t$ for $z$ from the Carath\'eodory kernel $\Omega(t_0)$ of
$\Omega(t)$, and derived a generalization of the L\"owner
equation. Ch.~Pommerenke \cite{Pom1} proposed to consider
subordination chains of domains that led him to a general
evolution equation. We mention here also papers by
V.~Gutlyanski{\u\i} \cite{Gut} and V.~Goryainov \cite{Gor} in this
direction. One can learn more about this method in monographs
\cite{Alexandrov, Duren, Pom2} (see also the references therein).
Let us draw reader's attention to Goryainov's approach \cite{Gor}
who suggested to use a method of semigroups to derive several
other parametric representations of classes of analytic maps and
to apply them to study  dynamics of stochastic branching
processes. This approach is based on the study of one-parameter
semi-flows on semigroups of conformal maps, their infinitesimal
descriptions, and evolution equations (see also \cite{Shoikhet}).

In 1959 Shah Dao-Shing \cite{Shah} suggested a parametric method
for quasiconformal automorphisms of $U$.  In another form this
method appeared in the paper by F.~Gehring and E.~Reich
\cite{Gehring}, and then, in \cite{Lawryn}. Later, Cheng Qi He
\cite{He} obtained an analogous equation for classes of
quasiconformally extendable univalent functions  (to be more
precise, in terms of inverse functions). Unlike the parametric
method for conformal maps, its analogue for quasiconformal maps
did not receive so much attention.

Several attempts have been launched to specialize the
L\"owner-Kufarev equation to obtain conformal maps that admit
quasiconformal extensions (see \cite{Becker, Becker1, Becker2,
Gut2}).
\s

Surprisingly, an analogous equation appeared in Fluid Dynamics in
the study of plane free boundary problems, where the time
dependence of the phase domain $\Omega(t)$ in a Hele-Shaw cell was
described by a one-parameter chain of univalent maps satisfying an
equation that now is known as the Polubarinova-Galin equation. It
appeared in the pioneering works by P.~Ya.~Polubarinova-Kochina
\cite{Polub1, Polub2} and L.~A.~Galin \cite{Galin} (see surveys
\cite{How}, \cite{Vas}). In contrast to the classical
L\"owner-Kufarev equation the latter is a non-linear (even
non-quasilinear) integro-differential equation and many elegant
properties of the L\"owner-Kufarev equation are less clear for
Polubarinova-Galin's one. A typical feature of the Hele-Shaw flow
is that starting with a  simply connected phase domain $\Omega(0)$
with a smooth boundary possible cusps may be developed during time
evolution. They are caused by  vanishing boundary derivatives as
well as by  topology change.

The principal goal of our paper is to study evolution equations
for conformal maps  with quasiconformal extensions. In particular,
we are interested in maps smoothly extendable onto the unit
circle. Our approach is based on the study of  evolutions on the
universal Teichm\"uller space $T$ and on the manifold $\Diff
S^1/\Rot S^1$ embedded into $T$.  Another question we are
interested in is how a Hele-Shaw flow is seen on the universal
Teichm\"uller space as on a general parametric space.

\section{The L\"owner-Kufarev and Polubarinova-Galin equations}

Let us consider a {\it subordination chain}  of simply connected
hyperbolic domains $\Omega(t)$ in the Riemann sphere $\hat{\mathbb
C}$, which is defined for $0\leq t< t_0$. This means that
$\Omega(t)\subset \Omega(s)$ whenever $t<s$. We suppose that all
$\Omega(t)$ are unbounded and $\infty\in \Omega(t)$ for all $t$.
By the Riemann Mapping Theorem we construct a subordination chain
of mappings $f(\zeta,t)$, $\zeta\in U^*$,  where each function
$\displaystyle
f(\zeta,t)=\alpha(t)\zeta+a_0(t)+\frac{a_1(t)}{\zeta}+\dots$ is a
meromorphic univalent map of $U^*$ onto $\Omega(t)$ for every
fixed $t$. Ch.~Pommerenke \cite{Pom1, Pom2} first introduced such
chains in order to generalize L\"owner's equation. His result says
that given a subordination chain of domains $\Omega(t)$ defined
for $t\in [0,t_0)$ with a differentiable real-valued coefficient
$\alpha(t)$ (in particular, $e^{-t}$ ), there exists  an analytic
regular function
$$p(\zeta,t)=p_0(t)+\frac{p_1(t)}{\zeta}+\frac{p_2(t)}{\zeta^2}+\dots,\quad
\zeta\in U^*,$$ such that $\R p(\zeta, t)>0$ in $\zeta\in U^*$ and
\begin{equation}
\frac{\partial f(\zeta,t)}{\partial t}=-\zeta\frac{\partial
f(\zeta,t)}{\partial \zeta}p(\zeta,t),\label{LKP}
\end{equation}
for almost all $t\in [0,t_0)$. The coefficient
$\alpha(t)=\alpha(0)\exp(-\int_0^tp_0(\tau)d\tau)$ is the
conformal radius of $\Omega(t)$. This equation now-a-days is known
as the L\"owner-Kufarev equation due to the contribution by
K.~L\"owner \cite{Loewner} and P.~P.~Kufarev \cite{Kufarev}.

\s
\noindent We consider two main {\bf questions:}
\begin{itemize}
\item If $\partial \Omega(t)$ is a quasicircle, what is $p(\zeta,
t)$?

\item If $\partial \Omega(t)$ is a smooth Jordan curve, what is $p(\zeta,
t)$?
\end{itemize}

\s

We draw reader's attention to the case of smooth boundaries and
their connection to free boundary problems of fluid dynamics. In
1898 H.~S.~Hele-Shaw \cite{Hele} proposed his famous cell that was
a device for investigating a flow of viscous fluid in a narrow gap
between two parallel plates.

The  dimensionless model of a moving viscous incompressible fluid
in the Hele-Shaw cell is described by a potential flow with the
velocity field $\Vb=(V_1, V_2)$. The pressure $p$ is the potential
for the fluid velocity $$\Vb=-\frac{h^2}{12\mu}\nabla p,$$ where
$h$ is the cell gap and $\mu$ is the viscosity of the fluid (see,
e.g. \cite{Ock, Saffman}).
 Through the similarity in the governing
equations, Hele-Shaw flows can be used to study  models of
saturated flows in porous media governed by Darcy's law. Over the
years various particular cases of such a flow have been
considered. Different driving mechanisms were employed, such as
surface tension or external forces (suction, injection). We
mention here a 600-paper bibliography of free and moving boundary
problems for Hele-Shaw and Stokes flows since 1898 up to 1998
collected by K.~A.~Gillow and S.~D.~Howison \cite{Gill}.

Since the work by Hele-Shaw several principal steps have been
made. Among them we distinguish the papers by
P.~Ya.~Polubarinova-Kochina \cite{Polub1, Polub2} and L.~A.~Galin
\cite{Galin} who suggested in 1945 a complex variable approach
that now is one of the basic tools for investigating the Hele-Shaw
evolution.

Let us consider the flow of a viscous fluid in a plane Hele-Shaw
cell under injection through a unique well which is placed at
infinity. Suppose that at the initial moment the phase domain
$\Omega_0$ occupied by the fluid is simply connected and bounded
by a smooth analytic curve $\Gamma_0$. This model can be thought
of as a receding air bubble in a viscous flow. The evolution of
the phase domains $\Omega(t)$ is described by an auxiliary
conformal mapping $f(\zeta,t)$ of $U^*$ onto $\Omega(t)$,
$\Omega(0)=\Omega_0$, normalized by
$f(\zeta,t)=\alpha(t)\zeta+a_0(t)+\frac{a_1(t)}{\zeta}+\dots,$
$\alpha(t)>0$. Here we denote the derivatives by $f'=\partial
f/\partial\zeta$, $\dot{f}=\partial f/\partial t$, and $t$ is the
time parameter. This mapping satisfies the equation
\begin{equation}
\R \Big[\dot f(\zeta,t)\overline{\zeta f'(\zeta,t)}\Big] =
-1,\quad \zeta=e^{i\theta}, \label{PG}
\end{equation}
under  suitable rescaling. L.~A.~Galin \cite{Galin},
P.~Ya.~Polubarinova-Kochina \cite{Polub1, Polub2} first derived
the equation (\ref{PG}) and stimulated deep investigation in
complex variable approach to free boundary problems (see, e.g.,
\cite{How, Vas} and the references therein).

 From (\ref{PG}) one can derive a L\"owner-Kufarev type
equation by the Schwarz-Poisson formula:
\begin{equation}
 \dot{f}=-\zeta f'\frac{1}{2\pi}\int\limits_{0}^{2\pi}\frac{1}{|f'(e^{i\theta},t)|^2}
\frac{\zeta+e^{i\theta}}{\zeta-e^{i\theta}}d\theta,\label{PGL}
\end{equation}
where $\zeta\in U^*$.

The equation (\ref{PG}) is equivalent to the kinematic condition
on the free boundary and, in particular, implies that the phase
domains $\Omega(t)$ form a subordination chain.  In contrary to
the classical L\"owner-Kufarev equation (\ref{LKP}) the equation
(\ref{PGL}) even is not quasilinear and the problem of the
short-time existence and uniqueness of the solution is much more
difficult. First it was proved by Yu.~P.~Vinogradov, P.~P.~Kufarev
\cite{VinKuf} in 1948 and later in 1993 by M.~Reissig, L.~Von
Wolfersdorf \cite{Reissig}. In fact, starting with a smooth domain
$\Omega_0$ the solution to (\ref{PG}) exists and unique locally in
time. It is known that the domains $\Omega(t)$ remain to have
smooth (even analytic) boundaries up to the time $t_0$ when
possible cusps are developed or the domain fails to be simply
connected. This means, in particular, that $\Omega(t)$ fails to be
a quasidisk as $t\to t_0$.

\s \noindent We ask the following {\bf question}: given an
 initial smooth phase domain
$\Omega_0$, and the Hele-Shaw evolution $\Omega(t)$, what kind of
evolution it produces on the universal Teichm\"uller space as a
natural general parametric space?

\begin{figure}
\centering
{\scalebox{0.8}{\includegraphics{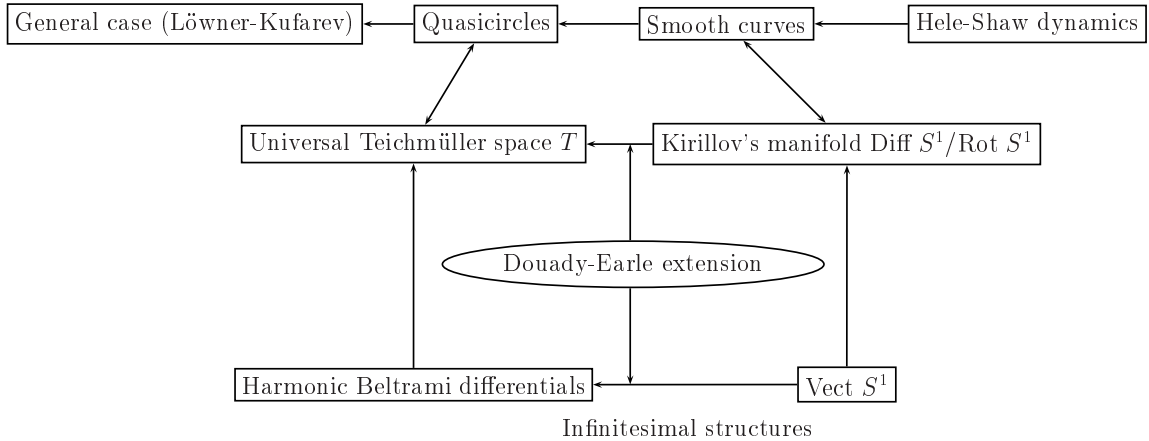}}}\caption[]{A general
scheme of study}\label{fig2}
\end{figure}

A general scheme of the proposed study is shown in
Figure~\ref{fig2}

\section{Infinitesimal structures of the universal Teichm\"uller space $T$}

Let us consider the family $\mathcal F$ of all quasiconformal
automorphisms of $U$. Every such  map $f$ satisfies the Beltrami
equation $f_{\bar{\zeta}}=\mu_f(\zeta)f_{\zeta}$ in $U$ in the
distributional sense, where $\mu_f$ is a measurable essentially
bounded function ($L^{\infty}(U)$) in $U$,
$\|\mu_f\|=\ess\sup_{U}|\mu_f(\zeta)|_{\infty}<1$. Conversely, for
each measurable Beltrami coefficient $\mu$ essentially bounded as
above, there exists a quasiconformal automorphism of $U$, that
satisfies the Beltrami equation, which is unique if  provided with
some conformal normalization, e.g., three point normalization $f(\pm
1)=\pm 1$, $f(i)=i$. Two normalized maps $f_1$ and $f_2$ are said to
be equivalent, $f_1\sim f_2$, if being extended onto the unit circle
$S^1$, the superposition $f_1\circ f_2^{-1}$ restricted to $S^1$ is
the identity map. The quotient set $\mathcal F/\sim$ is called {\it
the universal Teichm\"uller space} $T$\index{Teichm\"uller space}.
It is a covering space for all Teichm\"uller spaces of analytically
finite Riemann surfaces. By definition we have two realizations of
$T$: as a set of equivalence classes of quasiconformal maps and, due
to the relation between $\mathcal F/\sim$ and the unit ball
$B\subset L^{\infty}(U)$, as a set of equivalence classes of
corresponding Beltrami coefficients.

The normalized maps from $\mathcal F$ form a group $\mathcal F_0$
with respect to superposition and the maps that act identically on
$S^1$ form a normal subgroup $\mathcal I$. Thus, $T$ is the quotient
of $T=\mathcal F_0/\mathcal I$.

If $g\in \mathcal F$, $f\in \mathcal F_0$, then there exists a
M\"obius transformation $h$, such that $h\circ f\circ g^{-1}\in
\mathcal F_0$. Let us denote by $[f]\in T$ the equivalence class
represented by $f\in \mathcal F_0$. Then, one defines {\it the
universal modular group $\mathcal M$}, $\omega \in\mathcal M$,
$\omega:\,T\to T$, by the formula $\omega([f])=[h\circ f\circ
g^{-1}]$. Its subgroup  $\mathcal M_0$ of right translations on $T$
is defined by $\omega_0([f])=[f\circ g^{-1}]$, where  $f,g\in
\mathcal F_0$.

An important fact (see \cite[Chapter III, Theorem 1.1]{Lehto}) is
that there are real analytic mappings in any equivalence class
$[f]\in T$.

Given a Beltrami coefficient $\mu\in B\subset L^{\infty}(U)$
 let us extend it by zero into $U^*$. We normalize the corresponding quasiconformal map $f$, which
is conformal in $U^*$, by $f(\zeta)=\zeta+a_1/\zeta+\dots$ about
infinity. Then, two Beltrami coefficients $\mu$ and $\nu$ are
equivalent if and only if the corresponding normalized mappings
$f^{\mu}$ and $f^{\nu}$ map $U^*$ onto one and the same domain in
$\overline{\mathbb C}$. Thus, the universal Teichm\"uller space can
be thought of as the family of all normalized conformal maps of
$U^*$ admitting quasiconformal extension. Moreover, any compact
subset of $T$ consists of conformal maps $f$ of $U^*$ that admit
quasiconformal extension to $U$ with $\|\mu_f\|_{\infty}\leq k<1$
for some $k$.

As we mentioned above, a normalized conformal map $f\in [f]\in T$
defined in $U^*$ can have a quasiconformal extension to $U$ which is
real analytic in $U$, but on the unit circle $f$ may behave quite
irregularly. For example, the resulting quasicircle $f(S^1)$ can
have the Hausdorff dimension greater than 1. \s

\noindent {\bf Remark.}
 Given a bounded $K$-quasicircle $\Gamma$, $K=(1+k)/(1-k)$, in the
plane let $N(\varepsilon,\Gamma)$ denote the minimal number of disks
of radius $\varepsilon>0$ that are needed to cover $\Gamma$. Let
$$\beta(K)={\rm sup}\sb \Gamma\,{\rm lim sup}\sb
{\varepsilon\rightarrow0}\,{\rm log} N(\varepsilon,\Gamma)/{\rm
log}(1/\varepsilon)$$ denote the supremum of the Minkowski dimension
of curves $\Gamma$ where $\Gamma$ ranges over all bounded
$K$-quasicircles. The Hausdorff dimension of $\Gamma$ is bounded
from above by $\beta(K)$ (see \cite{BeckPom}). In \cite{BeckPom} it
was also established several explicit estimates for $\beta(K)$ ,
e.g., $\beta(K)\leq2-cK\sp {-3.41}$. \s

Let us denote by ${\Sigma}_0^{qc}\subset \Sigma_0$ the class of
those univalent conformal maps $f$ defined in $U^*$ which admit a
 quasiconformal extension to $U$, normalized by
$f(\zeta)=\zeta+a_1/\zeta+\dots$. Let $x,y\in T$ and $f,g\in
{\Sigma}_0^{qc}$ be such that $\mu_f\in x$ and $\mu_g\in y$. Then,
{\it the Teichm\"uller distance} $\tau(x,y)$ on $T$ is defined as
$$\tau(x,y)=\inf_{\mu_f\in x,\,\,\mu_g\in
y}\frac{1}{2}\log\frac{1+\|\mu_{g\circ
f^{-1}}\|_{\infty}}{1-\|\mu_{g\circ f^{-1}}\|_{\infty}}.$$

For a given $x\in T$ we consider an extremal Beltrami coefficient
$\mu^*$ such that $\|\mu^*\|_{\infty}=\inf_{\nu\in
x}\|\nu\|_{\infty}$. Let us remark that $\mu^*$ need not be unique.
A geodesic on $T$ can be described in terms of the extremal
coefficient $\mu^*$ as a continuous homomorphism $x_t:[0,1]\mapsto
T$ such that $\tau(0,x_t)=t\tau(0,x_1)$. Due to the above remark the
geodesic need not be unique as well.

We consider the Banach space $B(U)$ of all functions holomorphic in
$U$ equipped with the norm $$\|\varphi\|_{B(U)}=\sup_{\zeta\in
U}|\varphi(\zeta)|(1-|\zeta|^2)^2.$$  For a function $f$ in $\Sigma$
the Schwarzian derivative $$ S_f(\zeta)=\frac{\partial}{\partial
\zeta}\left(\frac{f''(\zeta)}{f'(\zeta)}\right)-
\frac{1}{2}\left(\frac{f''(\zeta)}{f'(\zeta)}\right)^2$$ is defined
and Nehari's \cite{Nehari} estimate $\|S_f(1/\zeta)\|_{B(U)}\leq 6$
holds. Given $x\in T$, $\mu\in x$ we construct the mapping
$f^{\mu}\in \Sigma_0^{qc}$ and have the homeomorphic embedding $T\to
B(U)$ by the Schwarzian derivative.

The universal Teichm\"uller space $T$ is an analytic infinite
dimensional Banach manifold modelled on $B(U)$. The Banach space
$B(U)$ is an infinite dimensional vector space that can be thought
of as the cotangent space to $T$ at the initial point (represented
by $\mu\equiv 0$). More rigorously, let the map $f^{\mu}$ be a
quasiconformal homeomorphism of the unit disk $U$. It has a
Fr\'echet derivative with respect to $\mu$ in a direction $\nu$. Let
us construct the variation of $f^{\tau\nu}\in \Sigma_0^{qc}$,
$\mu=\tau\nu$, with respect to a small parameter $\tau$:
\[
f^{\tau\nu}(\zeta)=\zeta+\tau V(\zeta)+o(\tau),\quad \zeta\in U^*.
\]
Taking the Schwarzian derivative in $U^*$ we get
\[
S_{f^{\tau\nu}}=\tau V'''(\zeta)+o(\tau),\quad \zeta\in U^*,
\]
locally uniformly in $U^*$. Taking into account the normalization of
the class $\Sigma_0^{qc}$ we have (see, e.g., \cite{Lehto})
$$V(\zeta)=-\frac{1}{\pi}\iint\limits_{U}\frac{\nu(w)d\sigma_{w}}{w-\zeta},\quad
V'''(\zeta)=-\frac{6}{\pi}\iint\limits_{U}\frac{\nu(w)d\sigma_{w}}{(w-\zeta)^4}.
$$
The integral formula implies $V'''(A(\zeta))A'(\zeta)^2=V'''(\zeta)$
(subject to the relation for the Beltrami coefficient
$\nu(A(\zeta))\overline{A'(\zeta)}=\nu(\zeta)A'(\zeta)$) for any
M\"obius transform $A$. Now let us change variables $\zeta\to
1/\bar{\zeta}$ and reduce the first variation to a holomorphic
function in the unit disk by changing $f^{\tau\nu}(\zeta)$ to
$g^{\tau\nu}(\zeta)\equiv\overline{{f}^{\tau\nu}(1/\bar{\zeta})}$.
Setting $\Lambda_{\nu}(\zeta)=S_{g^{\tau\nu}}(\zeta)$ and
$\dot{\Lambda}_{\nu}(\zeta)=\frac{1}{\zeta^4}\overline{V'''(1/\bar{\zeta})}$
we have (see, e.g., \cite[Section 6.5, Theorem 5]{GL}) that
\[
\Lambda_{\nu}(\zeta)-\tau
\dot{\Lambda}_{\nu}(\zeta)=\frac{o(\tau)}{(1-|\zeta|^2)^2}.
\]
So the operator $\dot{\Lambda}_{\nu}$ is the derivative of
$\Lambda_{\nu}$ at the initial point of the universal Teichm\"uller
space with respect to the norm of the Banach space $B(U)$.   The
reproducing property of the Bergman integral gives
\begin{equation}
\varphi(\zeta)=\frac{3}{\pi}\iint\limits_{U}\frac{\varphi(w)
(1-|w|^2)^2d\sigma_{w}}{(1-\bar{w}\zeta)^4},\quad\varphi\in
B(U).\label{aaa}
\end{equation}
The latter integral leads us to the so-called harmonic (Bers')
Beltrami differential
\[
\nu(\zeta)=\Lambda^*_{\varphi}(\zeta)\equiv-\frac{1}{2}
\overline{\varphi(\zeta)}(1-|\zeta|^2)^2,\quad \zeta\in U.
\]
Let us denote by $A(U)$ the Banach space of analytic functions with
the finite $L^1$ norm in the unit disk. We have that
$A(U)\hookrightarrow B(U)$ is a continuous inclusion (see, e.g.,
\cite[Section 1.4.2]{Nagbook}). On $L^{\infty}(U)\times A(U)$ one
can define a coupling
$$ \left<\mu,\varphi\right>:= \iint\limits_{U}\mu(\zeta)
\varphi(\zeta)\, d\sigma_{\zeta}. $$  Denote by $N$ the space of
{\it locally trivial Beltrami coefficients}\index{locally trivial
Beltrami coefficient}, which is the subspace of $L^{\infty}(U)$ that
annihilates  the operator $ \left<\cdot,\varphi\right>$ for all
$\varphi\in A(U)$. Then, one can identify the tangent space to $T$
at the initial point with the space $H:= L^{\infty}(U)/N$. It is
natural to relate it to a subspace of $L^{\infty}(U)$. The
superposition $\dot{\Lambda}_{\nu}\circ \Lambda^*_{\varphi}$ acts
identically on $A(U)$ due to (\ref{aaa}).  The space $N$ is also the
kernel of the operator $\dot{\Lambda}_{\nu}$.   Thus, the operator
$\Lambda\sp*$ splits the following exact sequence $$ 0
\longrightarrow N\hookrightarrow
L^{\infty}(U)\stackrel{\dot{\Lambda}_{\nu}}{\longrightarrow}
A(U)\longrightarrow 0. $$ Then, $H= \Lambda\sp*(A(U))\cong
L^{\infty}(U)/N$.  The coupling $\langle\mu, \varphi\rangle$ defines
$A(U)$ as a cotangent space. Let $A^2(U)$ denote the Banach space of
analytic functions $\varphi$ with the finite norm
\[
\|\varphi\|_{A^2(U)}=\iint\limits_{U}
|\varphi(\zeta)|^2(1-|\zeta|\sp2)\sp2 d\sigma_{\zeta}.
\]
Then $A(U)\hookrightarrow A^2(U)$ and  Petersson's Hermitian product
\cite{Petersson} is defined on $A^2(U)$ as
$$ (\varphi_1,\varphi_2)= \iint\limits_{U}
\varphi_1(\zeta)\overline{\varphi_2(\zeta)}(1-|\zeta|\sp2)\sp2
d\sigma_{\zeta}. $$ The K\"ahlerian Weil-Petersson metric
$\{\nu_1,\nu_2\}=\langle\nu_1,\dot{\Lambda}_{\nu_2}\rangle$ can be
defined on the tangent space to $T$ and gives a K\"ahlerian manifold
structure to $T$.

The universal Teichm\"uller space is a smooth manifold on which a
Lie group $\Diff T$ of real sense preserving diffeomorphisms  is
defined. The tangent bundle is defined on $T$ and is represented by
the harmonic differentials from $H$ translated to all points of $T$.
We will consider tangent vectors from $H$ at the initial point of
$T$ represented by the map $f(\zeta)\equiv \zeta$. The
Weil-Petersson metric defines a Lie algebra  of vector fields on $T$
by the Poisson-Lie bracket
$[\nu_1,\nu_2]=\{\nu_2,\nu_1\}-\{\nu_1,\nu_2\}$, where
$\nu_1,\nu_2\in H$. One can define the Poisson-Lie bracket at all
other points of $T$ by left translations from $\Diff T$. To each
element $[x]$ from $\Diff T$ an element $x$ from $T$ is associated
as an image of the initial point. Therefore, a curve in $\Diff T$
generates a traced curve in $T$ that can be realized by a
one-parameter family of quasiconfromal maps from $\Sigma_0^{qc}$.

For each tangent vector $\nu\in H$ there is a one-parameter
semi-flow in $\Diff T$ and a corresponding flow $x^{\tau}\in T$ with
the velocity vector $\nu$. To make an explicit representation we use
the variational formula for the subclass $\Sigma_0^{qc}$ of
$\Sigma_0$ of functions with quasiconformal extension (see, e.g.,
\cite{Lehto}) to $\overline{\mathbb C}$. If $f^{\mu}\in
\Sigma_0^{qc}$, $\nu\in H$ and
$$\mu_f(\zeta,\tau)=\left\{
\begin{array}{ll}
\tau\nu(\zeta)+o(\tau) &\mbox{if $\zeta\in U$,}\\ 0 & \mbox{if
$\zeta\in U^*$,}
\end{array}\right.
$$ then the map $$f^{\mu}(\zeta)=\zeta-\frac{\tau
}{\pi}\iint\limits_{U}\frac{\nu(w)d\sigma_{w}}{w-\zeta}+o(\tau)
$$ locally describes the semi-flow $x^{\tau}$ on $T$.

\section{$\Diff S^1/\Rot S^1$ embedded into $T$}

In this section we study a diffeomorphic embedding of the
homogeneous manifold $\Diff S^1/\Rot S^1$ into the universal
Teichm\"uller space $T$.

\subsection{Homogeneous manifold $\Diff S^1/\Rot S^1$}

We denote the Lie group of $C^{\infty}$ sense preserving
diffeomorphisms of the unit circle $S^1$  by $\Diff S^1$\index{Diff
$S^1$/Rot $S^1$}. Each element of $\Diff S^1$ is represented as
$z=e^{i\phi(\theta)}$ with a monotone increasing, $C^{\infty}$
real-valued function $\phi(\theta)$, such that
$\phi(\theta+2\pi)=\phi(\theta)+2\pi$.
 The Lie algebra for $\Diff S^1$ is
identified with the Lie algebra $\Vect S^1$ of smooth ($C^{\infty}$)
tangent vector fields to $S^1$ with the Poisson - Lie
bracket\index{Poisson - Lie bracket} given by
$$[\phi_1,\phi_2]={\phi}_1{\phi}'_2-{\phi}_2{\phi}'_1. $$
Fixing  the trigonometric basis in $\Vect S^1$ the commutator
relations take the form
\begin{eqnarray*}
\left[\cos\,n\theta, \cos\,m\theta\right] & = &
\frac{n-m}{2}\sin\,(n+m)\theta+ \frac{n+m}{2}\sin\,(n-m)\theta,\\
\left[\sin\,n\theta, \sin\,m\theta\right]& = &
\frac{m-n}{2}\sin\,(n+m)\theta+ \frac{n+m}{2}\sin\,(n-m)\theta,\\
\left[\sin\,n\theta, \cos\,m\theta\right] & = &
\frac{m-n}{2}\cos\,(n+m)\theta- \frac{n+m}{2}\cos\,(n-m)\theta.
\end{eqnarray*}
There is no general theory of infinite dimensional Lie groups,
example of which is under consideration.  The interest to this
particular case comes first of all from the string  theory where the
Virasoro algebra\index{Virasoro algebra} appears as the central
extension of $\Vect S^1$. Entire necessary background for the
construction of the theory of unitary representations of $\Diff S^1$
is found in the study of Kirillov's homogeneous K\"ahlerian manifold
$M=\Diff S^1/\Rot S^1$, where $\Rot S^1$ denotes the group of
rotations of $S^1$. The group $\Diff S^1$ acts as a group of
translations on the manifold $M$ with $\Rot S^1$ as a stabilizer.
The K\"ahlerian geometry of $M$ has been described by Kirillov and
Yuriev in \cite{KY1}. The manifold $M$ admits several
representations, in particular, in the space of smooth probability
measures, symplectic realization in the space of quadratic
differentials. We will use its analytic representation that is based
on the class $\tilde{\Sigma}_0$ of functions from $\Sigma_0$ which
being extended onto the closure $\overline{U}^*$ of $U^*$  are
supposed to be smooth on $S^1$. The class $\tilde{\Sigma}_0$ is
dense in $\Sigma_0$ in the local uniform topology of $U^*$.

\begin{figure}
\centering
{\scalebox{0.8}{\includegraphics{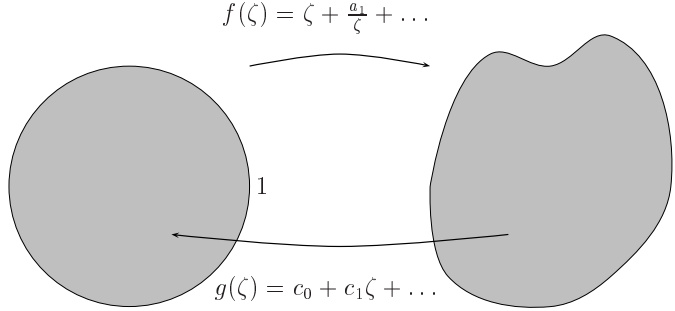}}}\caption[]{Representation
of $M$}\label{aafig3}
\end{figure}

Let $\tilde{S}$ denote the class of all univalent holomorphic maps
in the unit disk  $g(\zeta)=c_0+c_1\zeta+c_2\zeta^2+\dots$ which are
smooth on $S^1$. Then, for each $f\in \tilde{\Sigma}_0$ we have
$\infty\in f(U^*)$ and there is an {\it adjoint map} $g\in
\tilde{S}$ such that $\overline{\mathbb C}\setminus
f(U^*)=\overline{g(U)}$. The superposition $g^{-1}\circ f$
restricted to $S^1$ is in $M$ (see Figure \ref{aafig3}).
Reciprocally, for each element of $M$ there exist such $f$ and $g$.
A piece-wise smooth closed Jordan curve is a quasicircle if and only
if it has no cusps. So any function $f$ from $\tilde{\Sigma}_0$ has
a quasiconformal extension to $U$. By this realization the manifold
$M$ is naturally embedded into the universal Teichm\"uller space
$T$. Moreover, the K\"ahlerian structure on $M$ corresponds to the
K\"ahlerian structure\index{K\"ahlerian manifold} on $T$ given by
the Weil-Petersson metric\index{Weil-Petersson metric}.

The Goluzin-Schiffer variational formulae\index{Goluzin-Schiffer
variational formula} lift the actions from the Lie algebra $\Vect
S^1$ onto $\tilde{\Sigma}_0$. Let $f\in\tilde{\Sigma}_0$ and let
$d(e^{i\theta})$ be a $C^{\infty}$ real-valued function in
$\theta\in(0,2\pi]$ from $\Vect S^1$ making an infinitesimal action
as $\theta \mapsto \theta+\tau d(e^{i\theta})$. Let us consider a
variation of $f$ given by
\begin{equation}
\delta_{d}f(\zeta)=\frac{-1}{2\pi
i}\int\limits_{S^1}\left(\frac{wf'(w)}{f(w)}\right)^2\frac{wd
(w)dw}{f(w)-f(\zeta)} .\label{var}
\end{equation}
 Kirillov and Yuriev \cite{KY1}, \cite{KY2} have established
that the variations $\delta_{d}f(\zeta)$ are closed with respect to
the commutator and the induced Lie algebra is the same as $\Vect
S^1$. Moreover, Kirillov's result \cite{Kir} states that there is
the exponential map $\Vect S^1\to \Diff S^1$ such that the subgroup
$\Rot S^1$ coincides with the stabilizer of the map $f(\zeta)\equiv
\zeta$ from $\tilde{\Sigma}_0$.

\subsection{Douady-Earle extension}

Let $\varphi \colon S^1 \to S^1$ be a circle quasisymmetric
homeomorphism, i.e., a homeomorphism that possesses a quasiconformal
extension into $U$ (for a precise definition see, e.g.,
\cite{Lehto}). Then $\varphi$ has infinitely many quasiconformal
extensions into $U$, one of the most remarkable of which is the
Beurling-Ahlfors extension \cite{BlAhl}. In 1986 Douady and Earle
\cite{DE} defined for any such $\varphi \colon S^1 \to S^1$ a
conformally natural extension $h \colon \overline{U} \to
\overline{U}$ from $\mathcal F$. The map $h$ is a homeomorphism
which is real analytic in the interior. The idea was to introduce
the concept of a conformal barycenter of a measure on $S^1 =
\partial{U}$. Douady and Earle  proved that $w=h(\zeta)\in \mathcal F$
satisfies the functional equation
\begin{equation}
F(\zeta,w) \equiv {\frac{1}{2\pi}}\int\limits\sb {S^1}\left(
{\frac{\varphi(z)-w}{1-\overline{w}\varphi(z)}}\right)
{\frac{1-{|\zeta|}\sp 2}{{|\zeta-z|}\sp 2}} |d z| = 0.\label{DE}
\end{equation}
An advantage of this extension is that if $\sigma,
\tau\in$M\"ob$(U)$, then the extension of
$\sigma\circ\varphi\circ\tau$ is given by $\sigma\circ h \circ\tau$,
what is not true for the Beurling-Ahlfors extension. The three-point
boundary normalization of $\mathcal F_0$ can be always attained, and
thus, the Douady-Earle extension\index{Douady-Earle extension} is
compatible with the definition of the universal Teichm\"uller space.
Later, in 1988, another proof of Douady-Earle's result has appeared
in \cite{Lecko} where the authors worked with the inverse function.
The functional equation (\ref{DE}), in particular, implies that a
$C^{\infty}$ mapping $\varphi$ representing an element from the
manifold $M$ has a real analytic extension $h\in \mathcal F$ which
is $C^{\infty}$ on $S^1$.

Let $f\in\tilde{\Sigma}_0$ represent an element from $\varphi\in M$.
Let $g\in \tilde{S}$ be the adjoint map, $g^{-1}\circ
f\Big|_{S^1}=\varphi$. If $h$ is the Douady-Earle extension of
$\varphi$, then $g\circ h\big|_{S^1}\equiv f\big|_{S^1}$ and $g\circ
h$ is a quasiconformal extension of $f\in\tilde{\Sigma}_0$. Given
$\varphi\in M$ we construct the mapping $f^{\mu}$ that satisfies the
normalization of the class $\tilde{\Sigma}_0$ and whose Beltrami
coefficient is
\begin{equation}
\mu_f(\zeta)=\frac{\overline{F_{\zeta}}F_{\bar{w}}-F_{\bar{\zeta}}\overline{F_w}}
{\overline{F_{\bar{\zeta}}}F_{\bar{w}}-F_{\zeta}\overline{F_w}},\quad
w=h(\zeta),\,\,\,\zeta\in U,\label{mu}
\end{equation}
with $\mu_f(\zeta)=0$ for $\zeta\in U^*$. The equivalence class
$[f^{\mu}]$ is a point of the universal Teichm\"uller space $T$. So
the Douady-Earle extension defines an explicit embedding of $M$ into
$T$.

\subsection{Semi-flows on $T$ and $M$}

As it was mentioned, the Weil-Petersson metric defines a Lie algebra
of vector fields on $T$ by the Poisson bracket
$[\nu_1,\nu_2]=\{\nu_2,\nu_1\}-\{\nu_1,\nu_2\}$, where
$\nu_1,\nu_2\in H$. One can define the Poisson bracket at all other
points of $T$ by left translations of the universal modular group.

We proceed restricting ourselves to $M$ embedded into $T$. The
complex form of Green's formula implies that (\ref{var}) for
$f(\zeta)\equiv \zeta$ is equivalent to
\begin{equation}
\delta_{d}\,\zeta=\frac{-1}{\pi}\iint\limits_{U}\frac{\partial_{\bar{w}}
(wd (w))d\sigma_w}{w-\zeta},\label{var1}
\end{equation}
where the distributional derivative $\partial_{\bar{w}}d(w)$ is
given in the unit disk $U$,  $d(w)$ is a continuous extension of the
$C^{\infty}$ function $d(e^{i\theta})\in \Vect S^1$ into $U$ that
has $L^s(U)$ distributional derivatives in $U$, $s>2$, and
$d\sigma_w$ is the area element in $U$. Thus, one can extract  the
elements from $H$ that are of the form
$\nu(\zeta)=\zeta\partial_{\bar{\zeta}}d(\zeta)$, where
$\partial_{\bar{\zeta}}$ means $\partial/\partial\bar{\zeta}$.

We are going to deduce an exact form of $\nu$ using the Douady-Earle
extension. For this we start with the variation of the element
$$\varphi(e^{i\theta},\tau)=e^{i\theta}(1+\tau
id(e^{i\theta}))+o(\tau),\quad \varphi\in M,\,\,\,d\in \Vect S^1,$$
and $\tau$ is small. The Beltrami coefficient of the extended
quasiconformal map $h$ has  its variation as
$\mu_h(\zeta)=\tau\nu(\zeta)+o(\tau)$, where
\begin{equation}
\nu(\zeta)=\frac{\frac{\partial}{\partial\,\tau}\left(\overline{{F}_{\zeta}^{\tau}}F^{\tau}_{\bar{w}}-
F^{\tau}_{\bar{\zeta}}\overline{{F}^{\tau}_{w}}\right)}
{\overline{{F}^{\tau}_{\bar{\zeta}}}F^{\tau}_{\bar{w}}-F^{\tau}_{\zeta}\overline{{F}^{\tau}_{w}}}\Bigg|_{\tau=0,\,w=\zeta},
\zeta\in U,\label{mu1}
\end{equation}
where
\begin{equation} F^{\tau}(\zeta,w)
= {\frac{1}{2\pi}}\int\limits\sb {S^1}\left(
{\frac{\varphi(z,\tau)-w}{1-\overline{w}\varphi(z,\tau)}}\right)
{\frac{1-{|\zeta|}\sp 2}{{|\zeta-z|}\sp 2}} |d z| = 0.\label{DE1}
\end{equation}
Thus, $\nu(\zeta)$ depends only on $d(e^{i\theta})$. We will give
explicit formulae in the next section. They can be obtained
substituting $\varphi(e^{i\theta},0)=e^{i\theta}$, and taking into
account that $$\overline{{F}_{\zeta}^{\tau}}F^{\tau}_{\bar{w}}-
F^{\tau}_{\bar{\zeta}}\overline{{F}^{\tau}_{w}}\Big|_{\tau=0,\,w=\zeta}=0.$$
The Lie algebra $\Vect S^1$ is embedded into the Lie algebra of $H$
by (\ref{mu1}), (\ref{DE1}). Hence, a flow given on $M$
corresponding to a vector $d\in \Vect S^1$ is represented as a flow
on the universal Teichm\"uller space $T$ corresponding to the vector
$\nu\in H$ given by (\ref{mu1}).

\section{Infinitesimal descriptions of semi-flows}

First of all we give an explicit formula that connects the vectors
$d(e^{i\theta})$ from $\Vect S^1$  with corresponding tangent
vectors $\nu(\zeta)\in H$ to the universal Teichm\"uller space $T$
making use of the Douady-Earle extension\index{Douady-Earle
extension}. These vectors give the infinitesimal description of
semi-flows on $M$ and $T$ respectively.

\begin{theorem} Let $d(e^{i\theta})\in \Vect S^1$ be the infintesimal
description of a flow $\varphi$ in $M$. Then, the corresponding
infinitesmial description $\nu(\zeta)\in H$ of this flow embedded
into $T$ is given by the function
\begin{equation}
\nu(\zeta)=\frac{3}{2\pi}\int\limits_{0}^{2\pi}\left(\frac{1-|\zeta|^2}{(1-e^{i\theta}\bar{\zeta})^2}\right)^2
e^{2i\theta}d(e^{i\theta})d\theta.\label{nnn}
\end{equation}
\end{theorem}
\begin{proof}
Let $\varphi(\zeta,\tau)=e^{i(\theta+\tau d(e^{i\theta}))}$ and
$h(\zeta,\tau)$ be the Douady-Earle extension of $\varphi$ into the
unit disk $U$, $\zeta\in U$  by means of (\ref{DE1}). If $\tau=0$,
then $h(\zeta,0)\equiv\zeta$. We calculate
\begin{eqnarray*}
\partial_{\zeta}
F^{\tau}(\zeta,w) & = & {\frac{1}{2\pi}}\int\limits_{0}^{2\pi}\left(
{\frac{\varphi(e^{i\theta},\tau)-w}{1-\overline{w}\varphi(e^{i\theta},\tau)}}\right)
\frac{e^{i\theta}(\bar{\zeta}-e^{-i\theta})^2}{|\zeta-e^{i\theta}|^4}d\theta,\\
\partial_{\bar{\zeta}}
F^{\tau}(\zeta,w) & = & {\frac{1}{2\pi}}\int\limits_{0}^{2\pi}\left(
{\frac{\varphi(e^{i\theta},\tau)-w}{1-\overline{w}\varphi(e^{i\theta},\tau)}}\right)
\frac{e^{-i\theta}(\zeta-e^{i\theta})^2}{|\zeta-e^{i\theta}|^4}d\theta,\\
\partial_{w}
F^{\tau}(\zeta,w) & = & {\frac{1}{2\pi}}\int\limits_{0}^{2\pi}\left(
{\frac{-1}{1-\overline{w}\varphi(e^{i\theta},\tau)}}\right)
\frac{1-|\zeta|^2}{|\zeta-e^{i\theta}|^2}d\theta,\\
\partial_{\bar{w}}
F^{\tau}(\zeta,w) & = & {\frac{1}{2\pi}}\int\limits_{0}^{2\pi}\left(
{\frac{\varphi(e^{i\theta},\tau)(\varphi(e^{i\theta},\tau)-w)}{(1-\overline{w}\varphi(e^{i\theta},\tau))^2}}\right)
\frac{1-|\zeta|^2}{|\zeta-e^{i\theta}|^2}d\theta.
\end{eqnarray*}
Substituting $\tau=0$ and $w=\zeta$ we have
\begin{eqnarray*}
\partial_{\zeta}
F^{\tau}(\zeta,w)\bigg|_{\tau=0, w=\zeta} & = &
\frac{1}{1-|\zeta|^2},\\
\partial_{\bar{\zeta}}
F^{\tau}(\zeta,w)\bigg|_{\tau=0, w=\zeta} & = & 0,\\
\partial_{w}
F^{\tau}(\zeta,w)\bigg|_{\tau=0, w=\zeta} & = &
\frac{-1}{1-|\zeta|^2},\\
\partial_{\bar{w}}
F^{\tau}(\zeta,w)\bigg|_{\tau=0, w=\zeta} & = & 0.
\end{eqnarray*}

We will use the properties of the Douady-Earle extension. Let us fix
a point $\zeta_0\in U$ and choose two M\"obius transformations
$\sigma, \delta$ of $U$ such that $\delta(0)=\zeta_0$ and
$\sigma(0)=h(\zeta_0,\tau)$. We set $g=\sigma^{-1}\circ h \circ
\delta$. Then, $g(0,\tau)=0$, $\dot{g}(0,\tau)=0$ and
\begin{eqnarray*}
\partial_{\zeta}g(0,\tau) & = &
\partial_{\zeta}h(\zeta_0,\tau)\frac{\delta'(0)}{\sigma'(0)},\\
\partial_{\bar{\zeta}}g(0,\tau) & = &
\partial_{\bar{\zeta}}h(\zeta_0,\tau)\frac{\overline{\delta'(0)}}{\sigma'(0)}.
\end{eqnarray*}
So we see that
$$\frac{\partial_{\bar{\zeta}}h(\zeta_0,\tau)}{\partial_{\zeta}h(\zeta_0,\tau)}=
\frac{\partial_{\bar{\zeta}}g(0,\tau)}{\partial_{\zeta}g(0,\tau)}\frac{\delta'(0)}{\overline{\delta'(0)}}.$$
By the property of the Douady-Earle extension we have that the
function $g(\zeta,\tau)$, $\zeta\in U$ is the extension of
$g(e^{i\theta},\tau)$ by means of (\ref{DE1}). If $\tau=0$, then
$g(\zeta,0)\equiv\zeta$. Now we put
$\psi(e^{i\theta},\tau)=g(e^{i\theta},\tau)$ in (\ref{DE1}) and
calculate variations in $\tau$
\begin{eqnarray*}
\frac{\partial}{\partial\tau}\partial_{\bar{\zeta}}
F^{\tau}(\zeta,w)\bigg|_{\tau=0, w=\zeta=0}  & = &
{\frac{1}{2\pi}}\int\limits_{0}^{2\pi}\left(
\frac{e^{i\theta}-\zeta}{(1-\overline{\zeta}e^{i\theta})^3}\right)
\frac{1-|\zeta|^2}{|\zeta-e^{i\theta}|^2}\frac{\partial\psi(e^{i\theta},\tau)}{\partial\tau}\bigg|_{\tau=0,\zeta=0}
d\theta \\  & = &
{\frac{1}{2\pi}}\int\limits_{0}^{2\pi}e^{i\theta}\dot{\psi}(e^{i\theta},0)d\theta,\\
\frac{\partial}{\partial\tau}\partial_{\bar{w}}
F^{\tau}(\zeta,w)\bigg|_{\tau=0, w=\zeta=0}  & = &
{\frac{1}{2\pi}}\int\limits_{0}^{2\pi}\left(
\frac{2e^{i\theta}-\zeta-|\zeta|^2e^{i\theta}}{(1-\overline{\zeta}e^{i\theta})^3}\right)
\frac{1-|\zeta|^2}{|\zeta-e^{i\theta}|^2}\frac{\partial\psi(e^{i\theta},\tau)}{\partial\tau}\bigg|_{\tau=0,\zeta=0}
d\theta \\
 & = & {\frac{1}{2\pi}}\int\limits_{0}^{2\pi}2e^{i\theta}\dot{\psi}(e^{i\theta},0)d\theta.
\end{eqnarray*}
Then, we can obtain the explicit form of the variation of the
Beltrami coefficient  by (\ref{mu1}) as
\begin{equation}
\frac{\partial}{\partial
\tau}\frac{\partial_{\bar{\zeta}}g(0,\tau)}{\partial_{\zeta}g(0,\tau)}\bigg|_{\tau=0}=
\frac{3}{2\pi}\int\limits_{0}^{2\pi}e^{i\theta}\dot{\psi}(e^{i\theta},0)d\theta.\label{dot}
\end{equation}

 The M\"obius transformation $\delta$ does not depend on $\tau$
whereas $\sigma$ does. Explicitly, we put $$\sigma^{-1}\circ h\circ
\delta(\zeta)
=\frac{h(\delta(\zeta),\tau)-h(\zeta_0,\tau)}{1-h(\delta(\zeta),\tau)\overline{h(\zeta_0,\tau)}},
$$ where
$\delta(\zeta)=(\zeta+\zeta_0)(1+\zeta\bar{\zeta}_0)^{-1}$. We
denote $e^{i\alpha}=\delta(e^{i\theta})$. Therefore, denoting by
$$
e^{i\alpha}=\delta(e^{i\theta})=\frac{e^{i\theta}+\zeta_0}{1+\bar{\zeta}_0e^{i\theta}},
$$ we have
$$\dot{g}(e^{i\theta},0)=\frac{\dot{h}(e^{i\alpha},0)(1-|\zeta_0|^2)-\dot{h}(\zeta_0,0)(1-\bar{\zeta}_0e^{i\alpha})+
\overline{\dot{h}(\zeta_0,0)}e^{i\alpha}(e^{i\alpha}-\zeta_0)}{(1-\bar{\zeta}_0e^{i\alpha})^2}.
$$ Then,
$$e^{i\theta}d\theta=\frac{1-|\zeta_0|^2}{(1-e^{i\alpha}\bar{\zeta}_0)^2}e^{i\alpha}d\alpha,
$$ and changing variables in (\ref{dot}), we obtain $$
\frac{\partial}{\partial
\tau}\frac{\partial_{\bar{\zeta}}g(0,\tau)}{\partial_{\zeta}g(0,\tau)}\bigg|_{\tau=0}=
\frac{3}{2\pi}\int\limits_{0}^{2\pi}\left(\frac{1-|\zeta_0|^2}{(1-e^{i\alpha}\bar{\zeta}_0)^2}\right)^2
e^{2i\alpha}d(e^{i\alpha})d\alpha.$$ Taking into account that
$\delta'(0)=1$ we come to the statement of the theorem.
\end{proof}
\begin{corollary} If
$q=\max\limits_{\theta\in[0,2\pi]}|d(e^{i\theta})|$, then
$$|\nu(\zeta)|\leq 3\frac{1+|\zeta|^2}{1-|\zeta|^2}q.$$
\end{corollary}
\begin{proof}
The formula given in the preceding theorem implies $$
\nu(\zeta)=\frac{3}{2\pi}\int\limits_{0}^{2\pi}\frac{1-|\zeta|^2}{(1-e^{i\alpha}\bar{\zeta})^2}
e^{i\alpha}d(\delta(e^{i\theta}))e^{i\theta}d\theta.$$ Changing
variables $\alpha\to \theta$ we obtain
\begin{equation}
\nu(\zeta)=\frac{3}{2\pi}\int\limits_{0}^{2\pi}\frac{e^{i\theta}+\zeta}{1+e^{i\theta}\bar{\zeta}}\frac{(1+e^{i\theta}\bar{\zeta})^2}{1-|\zeta|^2}
e^{i\theta}d(\delta(e^{i\theta}))d\theta.\label{aa}
\end{equation}
 Next, we obviously
estimate $|\nu|$ as in the statement of the corollary. \end{proof}

As we see, the given estimate is good enough when $|\zeta|$ is not
close to 1. Let us now give an asymptotic estimate for
$|\nu(\zeta)|$ in the case $|\zeta|\sim 1$.

\begin{corollary} There exists a constant $M$ independent of
$\zeta$such that $$|\nu(\zeta)|\leq M\frac{1-|\zeta|^2}{|\zeta|^2}.
$$ In particular, $|\nu(\zeta)|=O(1-|\zeta|^2)$ as $|\zeta|\sim 1$.
\end{corollary}
\begin{proof}
We integrate by parts the right-hand side in the formula (\ref{nnn})
twice and come to the following expression
\begin{equation}
\nu(\zeta)=-\frac{(1-|\zeta|^2)}{4\pi\bar{\zeta}^2}\int\limits_{0}^{2\pi}\frac{1-|\zeta|^2}
{(1-e^{i\theta}\bar{\zeta})^2}\left(i\frac{\partial[e^{i\theta}
d(e^{i\theta})]}{\partial \theta}+\frac{\partial^2
[e^{i\theta}d(e^{i\theta})]}{\partial \theta^2}\right)d\theta.
\end{equation}
The absolute value of the above integral is bounded because of the
Poisson kernel in it and due to the smoothness of the function $d$.
\end{proof}

\section{Parametric representation of univalent maps with quasiconformal extensions}

\subsection{Semigroups of conformal maps}

The basic ideas that we use in this section come from Goryainov's
works \cite{Gor}, \cite{Gor1} and the monograph by Shoikhet
\cite{Shoikhet}.

 We consider the semigroup $\mathcal G$ of conformal
univalent maps from $U^*$ into itself with composition as the
semigroup\index{semigroup} operation. This makes $\mathcal G$ a
topological semigroup with respect to the topology of local uniform
convergence on $U^*$. We impose the natural normalization for such
conformal maps: $\displaystyle\Phi(\zeta)=\beta \zeta+b_0
+\frac{b_1}{\zeta} +\dots$, $\zeta\in U^*$, $\beta>0$. The unit of
the semigroup is the identity. Let us construct on $\mathcal G$ a
one-parameter semi-flow\index{semi-flow} $\Phi^{\tau}$, that is, a
continuous homomorphism from $\mathbb R\sp +$ into $\mathcal G$,
with the parameter ${\tau}\geq 0$. For any fixed ${\tau}\geq 0$ the
element $\Phi^{\tau}$ is from $\mathcal G$ and is represented by a
conformal map $\displaystyle \Phi(\zeta,{\tau})=\beta(\tau)
\zeta+b_0(\tau)+\frac{b_1(\tau)}{\zeta}+\dots$ from $U^*$ onto the
domain $\Phi(U^*,{\tau})\subset U^*$. The element $\Phi^{\tau}$
satisfies the following properties:
\begin{itemize}
\item $\Phi^0= id$;

\item $\Phi^{{\tau}+s}=\Phi(\Phi(\zeta,{\tau}),s)$, for $\tau,
s\geq 0$;

\item $\Phi(\zeta,{\tau})\to \zeta$ locally uniformly in $U^*$ as
$\tau\to 0$.

\end{itemize}
In particular, $\beta(0)=1$. This semi-flow is generated by a vector
field $v(\zeta)$ if for each $\zeta\in U^*$ the function
$w=\Phi(\zeta,\tau)$, $\tau\geq 0$ is a solution of an autonomous
differential equation $dw/d\tau =v(w)$ with the initial condition
$w|_{\tau=0}=\zeta$. The semi-flow can be extended to a symmetric
interval $(-t,t)$ by putting $\Phi^{-\tau}= \Phi^{-1}(\zeta,\tau)$.
Certainly, the latter function is defined on the set
$\Phi(U^*,\tau)$. Admitting this restriction for negative $\tau$ we
define a one-parameter family $\Phi^{\tau}$ for $\tau\in (-t,t)$.

For a semi-flow $\Phi^{\tau}$ on $\mathcal G$ there is an
infinitesimal generator at ${\tau}=0$ constructed by the following
procedure. Any element $\Phi^{\tau}$ is represented by a conformal
map $\Phi(\zeta,{\tau})$ that satisfies the Schwarz
Lemma\index{Schwarz Lemma} for the maps $U^*\to U^*$, and hence,
$$\R\frac{\zeta}{\Phi(\zeta,{\tau})}\leq
\Big|\frac{\zeta}{\Phi(\zeta,{\tau})}\Big|\leq 1,\quad \zeta\in
U^*,$$ where the equality sign is attained only for $\Phi^0=id\simeq
\Phi(\zeta,0)\equiv \zeta$. Therefore, the following limit exists
(see, e.g., \cite{Gor}, \cite{Gor1}, \cite{Shoikhet})
$$\lim\limits_{{\tau}\to
0}\R\frac{\zeta-\Phi(\zeta,{\tau})}{{\tau}\Phi(\zeta,{\tau})}=-\R\frac{\frac{\partial
\Phi(\zeta,{\tau})}{\partial {\tau}}\Big|_{{\tau}=0}}{\zeta}\leq
0,$$ and the representation $$\frac{\partial
\Phi(\zeta,{\tau})}{\partial {\tau}}\Big|_{{\tau}=0}=\zeta
p(\zeta)$$ holds, where $\displaystyle p(\zeta)=p_0+p_1/\zeta+\dots$
is an analytic function in $U^*$ with  positive real part, and
\begin{equation}
\frac{\partial\beta(\tau)}{\partial
\tau}\Big|_{\tau=0}=p_0.\label{aagenerator}
\end{equation}
In \cite{Gor2} it was shown that $\Phi^{\tau}$ is even $C^{\infty}$
with respect to $\tau$. The function $\zeta p(\zeta)$ is an
infinitesimal generator for $\Phi^{\tau}$ at ${\tau}=0$, and the
following variational formula holds
\begin{equation}
\Phi(\zeta,\tau)=\zeta+\tau\,\zeta p(\zeta)+o(\tau),\quad
\beta(\tau)=1+\tau p_0+o(\tau).\label{aa1}
\end{equation}
The convergence is thought of as local uniform.  We rewrite
(\ref{aa1}) as
\begin{equation}
\Phi(\zeta,\tau)=(1+\tau p_0)\zeta+\tau\,\zeta(
p(\zeta)-p_0)+o(\tau)=\beta(\tau)\zeta+\tau\,\zeta(
p(\zeta)-p_0)+o(\tau).\label{aa2}
\end{equation}

Now let us proceed with the semigroup $\mathcal G^{qc}\subset
\mathcal G$ of quasiconformal automorphisms of $\overline{\mathbb
C}$.
 A quasiconformal map $\Phi$
representing an element of $\mathcal G^{qc}$ satisfies the Beltrami
equation in $\overline{\mathbb C}$
$$\Phi_{\bar{\zeta}}=\mu_{\Phi}(\zeta)\Phi_{\zeta},$$ with the distributional
derivatives $\Phi_{\bar{\zeta}}$ and $\Phi_{\zeta}$, where
$\mu_{\Phi}(\zeta)$ is a measurable function vanishing in $U^*$ and
essentially bounded in $U$ by
$$\|\mu_{\Phi}\|=\ess\sup_{U}|\mu_{\Phi}(\zeta)|\leq k <1,$$ for some
$k$. If $k$ is sufficiently small, then the function
$\displaystyle\frac{\Phi-b_0}{\beta}$ satisfies the variational
formula (see, e.g., \cite{Lehto})
\begin{equation}
\frac{\Phi(\zeta)-b_0}{\beta}=\zeta-\frac{1}{\pi}\iint\limits_{U}\frac{\mu_{\Phi}(w)d\sigma_{w}}{w-\zeta}+o(k),
\label{aa3}
\end{equation}
where $d\sigma_{w}$ stands for the area element in the $w$-plane.

Now for each $\tau$ small and $\Phi^{\tau}\in \mathcal G^{qc}$ the
mapping
$h(\zeta,\tau)=\frac{\Phi(\zeta,\tau)-b_0(\tau)}{\beta(\tau)}$ is
from $\Sigma_0^{qc}$ and represents an equivalence class
$[h^{\tau}]\in T$. Consider the one-parameter curve $x^{\tau}\in T$
that corresponds to $[h^{\tau}]$ and a velocity vector
$\nu(\zeta)\in H$ (that is not trivial), such that
$$\mu_h(\zeta,\tau)=\mu_{\Phi}(\zeta,\tau)=\tau\nu(\zeta)+o(\tau).$$ We take
into account that $\Phi(\zeta,0)\equiv \zeta$ in $U^*$ and is
extended up to the identity map of $\overline{\mathbb C}$.

The formula (\ref{aa3}) can be rewritten for $\Phi(\zeta,\tau)$ as
\begin{equation}
\frac{\Phi(\zeta,\tau)-b_0(\tau)}{\beta(\tau)}=
\zeta-\frac{\tau}{\pi}\iint\limits_{U}\frac{\nu(w)d\sigma_{w}}{w-\zeta}+o(\tau).\label{aa4}
\end{equation}
Comparing with (\ref{aa2}) we come to the conclusion about $\Phi$:
\begin{equation}
\Phi(\zeta,\tau)=\beta(\tau)\zeta+\tau
p_1-\frac{\tau}{\pi}\iint\limits_{U}
\frac{\nu(w)d\sigma_{w}}{w-\zeta}+o(\tau). \label{aa5}
\end{equation}
The relations (\ref{aa1}, \ref{aa2}, \ref{aa5}) imply that
\begin{equation}
p(z)=p_0+\frac{p_1}{\zeta}-\frac{1}{\pi}\iint\limits_{U}\frac{\nu(w)d\sigma_{w}}{\zeta(w-\zeta)}.
\label{aa6}
\end{equation}
The constants $p_0,p_1$ and the function $\nu$ must be such that $\R
p(z)>0$ for all $z\in U^*$.

We summarize these observations in  the following theorem.

\begin{theorem}
Let $\Phi^{\tau}$ be a semi-flow in $\mathcal G^{qc}$. Then it is
generated by the vector field $v(\zeta)=\zeta p(\zeta)$,
\[
p(z)=p_0+\frac{p_1}{\zeta}-\frac{1}{\pi}\iint\limits_{U}\frac{\nu(w)d\sigma_{w}}{\zeta(w-\zeta)},
\]
where $\nu(\zeta)\in H$ is a harmonic Beltrami differential and the
holomorphic function $p(\zeta)$ has positive real part in $U^*$.
\end{theorem}

This theorem implies that at any point $\tau\geq 0$ we have
\[
\frac{\partial\Phi(\zeta,\tau)}{\partial
\tau}=\Phi(\zeta,\tau)p(\Phi(\zeta,\tau)).
\]
\subsection{Evolution families and differential equations}

A subset $\Phi^{t,s}$ of $\mathcal G$, $0\leq s\leq t$ is called an
{\it evolution family} in $\mathcal G$ if
\begin{itemize}
\item $\Phi^{t,t}= id$;

\item $\Phi^{t,s}=\Phi^{t,r}\circ \Phi^{r,s}$, for $0\leq s\leq
r\leq t$;

\item $\Phi^{t,s}\to id$ locally uniformly in $U^*$ as
$t,s\to\tau$.

\end{itemize}
In particular, if $\Phi^\tau$ is a one-parameter
semi-flow\index{semi-flow}, then $\Phi^{t-s}$ is an evolution
family\index{evolution family}. We consider a subordination chain of
mappings $f(\zeta,t)$, $\zeta\in U^*$, $t\in [0,t_0)$, where the
function $\displaystyle
f(\zeta,t)=\alpha(t)z+a_0(t)+a_1(t)/\zeta+\dots$ is a meromorphic
univalent map\index{meromorphic univalent map} $U^*\to
\overline{\mathbb C}$ for each fixed $t$ and $f(U^*,s)\subset
f(U^*,t)$ for $s<t$. Let us assume that this subordination
chain\index{subordination chain} exists for $t$ in an interval
$[0,t_0)$.

Let us pass to the semigroup $\mathcal G^{qc}$. So $\Phi^{t,s}$ now
has a quasiconformal extension to $U$ and being restricted to $U^*$
is from $\mathcal G$. Moreover, $\Phi^{t,s}\to id$ locally uniformly
in $\mathbb C$ as $t,s\to \tau$.

For each $t$ fixed in $[0,t_0)$ the map $f(\zeta,t)$ has a
quasiconformal extension into $U$ (that can be assumed even real
analytic).  An important presupposition is that $f(\zeta,t)$
generates a {\it nontrivial path} in the universal Teichm\"uller
space $T$\index{Teichm\"uller space}. This means that for any
$t_1,t_2\in [0,t_0)$, $t_1\neq t_2$, the mapping $f(\zeta,t_2)$,
$\zeta\in U^*$, can not be obtained from $f(\zeta,t_1)$ by a
M\"obius transform, or taking into account the normalization of $f$,
by multiplying by a constant. We construct the superposition
$f^{-1}(f(\zeta,s),t)$ for $t\in [0,t_0)$, $s \leq t$. Putting
$s=t-\tau$ we denote this mapping by $\Phi(\zeta,t,\tau)$.

Now we suppose the following conditions for  $f(\zeta,t)$.
\begin{itemize}

\item[(i)] \ \ The maps $f(\zeta,t)$ form a subordination chain in
$U^*$, $t\in [0,t_0)$.

\item[(ii)] \ \ \ The map $f(\zeta,t)$ is holomorphic in $U^*$,
$f(\zeta,t)=\alpha(t)\zeta+a_0(t)+a_1(t)/\zeta+\dots$, where
$\alpha(t)>0$ and differentiable with respect to $t$.

\item[(iii)] \ \ \  The map $f(\zeta,t)$ is a quasiconformal
homeomorphism of $\overline{\mathbb C}$.

\item[(iv)] \ \ \  The chain of maps $f(\zeta,t)$ is not trivial.

\item[(v)] \ \  The Beltrami coefficient $\mu_f(\zeta,t)$ of this
map is differentiable with respect to $t$ locally uniformly in $U$,
vanishes in some neighbourhood of $U^*$ (independently of $t$).
\end{itemize}

\begin{figure}[ht]
\centering
{\scalebox{0.8}{\includegraphics{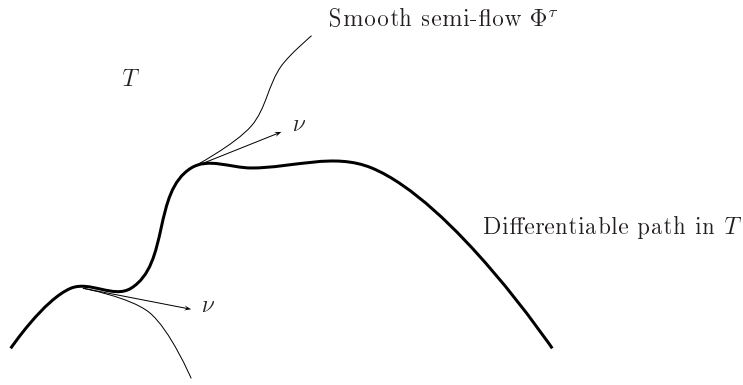}}}\caption[]{The pencil
of tangent smooth semi-flows}\label{aafig4}
\end{figure}

 The function
$\Phi(\zeta,t,\tau)$ is embedded into
 an evolution family in
$\mathcal G$. It is differentiable with regard to $\tau$ and $t$ in
$[0,t_0)$, and $\Phi(\zeta,t,0)=\zeta$. Fix $t$ and let $D_{\tau}=
\Phi^{-1}(U^*,t,\tau)\setminus U^*$. Then, there exists $\nu\in H$
such that the Beltrami coefficient $\mu$ is of the form
$\mu_{\Phi}(\zeta,t,\tau)=\tau\nu(\zeta,t)+o(\tau)$ in $U\setminus
D_{\tau}$, $\mu_{\Phi}(\zeta,t,\tau)=\mu_f(\zeta,t-\tau)$ in
$D_{\tau}$, and vanishes in $\hat{U^*}$. We make $\tau$ sufficiently
small such that $\mu_{\Phi}(\zeta,t,\tau)$ vanishes in $D_{\tau}$
too. Therefore, $\zeta=\lim_{\tau\to 0}\Phi(\zeta,t,\tau)$ locally
uniformly in $\mathbb C$ and $\Phi(\zeta,t,\tau)$ is embedded now
into an evolution family in $\mathcal G^{qc}$. The identity map is
embedded into a semi-flow $\Phi^{\tau}\subset {\mathcal G}^{qc}$
(which is smooth) as the initial point with the same velocity vector
 $$ \frac{\partial
\Phi(\zeta,t,{\tau})}{\partial {\tau}}\Big|_{{\tau}=0}=\zeta
p(\zeta,t),\quad \zeta\in U^*, $$ that leads to equation (\ref{LKP})
(the semi-flow $\Phi^{\tau}$ is tangent to the evolution family at
the origin). Actually, the differentiable trajectory $f(\zeta,t)$
generates a pencil of tangent smooth semi-flows with starting
tangent vectors $\zeta p(\zeta,t)$ (that may be only measurable with
respect to $t$). The projection to the universal Teichm\"uller space
is shown in Figure~\ref{aafig4}.

The requirement of non-triviality makes it possible to use the
variation (\ref{aa4}). Therefore, the conclusion is that the
function $f(\zeta,t)$ satisfies the equation (\ref{LKP}) where the
function $p(\zeta,t)$ is given by
\[
p(\zeta,t)=p_0(t)+\frac{p_1(t)}{\zeta}-\frac{1}{\pi}\iint\limits_{U}\frac{\nu(w,t)d\sigma_{w}}{\zeta(w-\zeta)},
\]
and has positive real part. The existence of $p_0(t), p_1(t)$ comes
from the existence of the subordination chain. We can assign the
normalization to $f(\zeta,t)$ controlling the change of the
conformal radius of the subordination chain by $e^{-t}$. In this
case, changing variables we obtain $p_0=1$, $p_1=0$.

Summarizing the conclusions about the function $p(\zeta,t)$ we come
to the following result.

\begin{theorem}\label{quasi} Let $f(\zeta,t)$ be a subordination  chain of
maps in $U^*$ that exists for $t\in [0,t_0)$  and satisfies the
conditions (i--v). Then, there are a real valued function
$p_0(t)>0$,  a complex valued function $p_1(t)$, and a harmonic
Beltrami differential $\nu(\zeta,t)$, such that $\R p(\zeta,t)>0$
for $\zeta\in U^*$,
$$p(\zeta,t)=p_0(t)+\frac{p_1(t)}{\zeta}-\frac{1}{\pi}\iint\limits_{U}\frac{\nu(w,t)d\sigma_{w}}{\zeta(w-\zeta)},
\quad \zeta\in U^*, $$ and $f(\zeta,t)$ satisfies the differential
equation
\begin{equation}
\frac{\partial f(\zeta,t)}{\partial t}=-\zeta \frac{\partial
f(\zeta,t)}{\partial \zeta}p(\zeta,t),\quad \zeta\in U^*,\label{aa7}
\end{equation}
in $t\in [0, t_0)$.
\end{theorem}

In the above theorem the function $\nu(\zeta,t)$ belongs to the
space of harmonic differentials. We ask now about another but
equivalent form of $\nu$  as well as whether one can extend the
equation $(\ref{aa7})$ onto the whole complex plane.

Writing $w=f(\zeta,t-\tau)$, $\Phi(\zeta,t,\tau)=f^{-1}(w,t)$ we
calculate the dilatation of the function $\Phi(\zeta,t,\tau)$ in
$U$. Note that $\Phi$ it is differentiable by $t,\tau$.
$$\mu_{\Phi}=\frac{\Phi_{\bar{\zeta}}}{\Phi_{\zeta}}=\frac{f^{-1}_{w}w_{\bar{\zeta}}+f^{-1}_{\bar{w}}\bar{w}_{\bar{\zeta}}}
{f^{-1}_{w}w_{\zeta}+f^{-1}_{\bar{w}}\bar{w}_{\zeta}}=
\frac{w_{\bar{\zeta}}+\mu_{f^{-1}}\bar{w}_{\bar{\zeta}}}
{w_{\zeta}+\mu_{f^{-1}}\bar{w}_{\zeta}}=
\frac{\bar{w}_{\bar{z}}}{w_{\zeta}} \frac{\mu_w\frac{\displaystyle
w_z}{\displaystyle \bar{w}_{\bar{\zeta}}}-\mu_{f}\frac{\displaystyle
f_{\zeta}}{\displaystyle \bar{f}_{\bar{\zeta}}}}
{1-\mu_{f}\overline{\mu_{w}}\frac{\displaystyle
f_{\zeta}\bar{w}_{\bar{\zeta}}}{\displaystyle
w_{\zeta}\bar{f}_{\bar{\zeta}}}}.$$ We use that $\mu_{f^{-1}}\circ
f=-\mu_f f_{\zeta}/\bar{f}_{\bar{\zeta}}$. Finally, $\mu_f$,
$f_{\zeta}$, $f_{\bar{\zeta}}$ are differentiable by $t$ almost
everywhere in $t\in [0,t_0)$, locally uniformly in $\zeta\in U$, and
$$\nu_0(\zeta,t)=\lim_{\tau\to
0}\frac{\mu_{\Phi}}{\tau}=-\frac{\bar{f}_{\bar{\zeta}}}{f_{\zeta}}\frac{\frac{\displaystyle\partial
}{\displaystyle\partial t}\left(\mu_f\frac{\displaystyle
f_{\zeta}}{\displaystyle
\bar{f}_{\bar{\zeta}}}\right)}{1-|\mu_f|^2},$$ where the limit
exists a.e. with respect to $t\in [0,t_0)$ locally uniformly in
$\zeta\in U$, or in terms of the inverse function
$$
\nu_0(\zeta,t)=\left(\frac{f^{-1}_w}{\bar{f}^{-1}_{\bar{w}}}\frac{\frac{\displaystyle\partial
\mu_{f^{-1}}}{\displaystyle\partial
t}}{1-|\mu_{f^{-1}}|^2}\right)\circ{f(\zeta,t)}.$$ Sometimes, it is
much better to operate just with dilatations, avoiding functions, so
we can rewrite the last expression as $$
\nu_0(z,t)=-\mu_f(z,t)\left[\frac{\frac{\displaystyle\partial
\log\mu_{f^{-1}}}{\displaystyle\partial
t}}{1-|\mu_{f^{-1}}|^2}\circ{f(z,t)}\right].$$ \s

\noindent{\it Remark.} The function $\nu(\zeta,t)$ in
Theorem~\ref{quasi} may be replaced by the function $\nu_0(\zeta,t)$
that belongs to the same equivalence class in $H$. \s

Let us consider one-parameter families of maps in $U^*$ normalized
by $\displaystyle
f(\zeta,t)=e^{-t}\zeta+\frac{a_1(t)}{\zeta}+\dots$. The inverse
result to the L\"owner-Kufarev equation\index{L\"owner-Kufarev
equation} states that given a holomorphic function
$p(\zeta,t)=1+p_1(t)/\zeta+\dots$ in $\zeta\in U^*$ with positive
real part the solution of the equation (\ref{aa7}) presents a
subordination chain (see, e.g., \cite{Pom2}). This enable us to give
a condition for $\nu_0$ that guarantees a normalized one-parameter
non-trivial family of maps $f(\zeta,t)$ to be a subordination chain

\begin{theorem} Let $f(\zeta,t)$ be a normalized one-parameter non-trivial family of maps
 for $\zeta\in U^*$ which satisfies
the conditions (ii--v) and is defined in an interval $[0,t_0)$. Let
each $f(\zeta,t)$ be a homeomorphism of $\overline{\mathbb C}$ which
is meromorphic in $U^*$, is normalized by $\displaystyle
f(\zeta,t)=e^{-t}\zeta+\frac{a_1(t)}{\zeta}+\dots$, and satisfies
(\ref{aa7}). Let the quasiconformal extension to $U$ be given by a
Beltrami coefficient $\mu_f=\mu(\zeta,t)$ which is differentiable
with respect to $t$ almost everywhere in $t\in[0,t_0)$. If
\[
\|\nu_0\|_{\infty}<\frac{\pi}{4\int_0^1s\K(s)ds}\approx
0.706859\dots,
\]
where ${\nu_0}(\zeta,t)$ is as  above  and $\K(\cdot)$ is the
complete elliptic integral, then $f(\zeta,t)$ is a normalized
subordination chain.
\end{theorem}
\begin{proof}
 Let $|\zeta|=\rho$, $w=re^{i\theta}$. We calculate
\begin{eqnarray*}
\Bigg|\frac{1}{\pi}\iint\limits_{U}\frac{\nu_0(w,t)d\sigma_w}{\zeta(w-\zeta)}\Bigg|
& \leq
&\frac{\|\nu_0\|_{\infty}}{\rho\pi}\iint\limits_{U}\frac{d\sigma_w}{|w-z|}
=
\frac{\|\nu_0\|_{\infty}}{\rho^2\pi}\iint\limits_{U}\frac{d\sigma_w}{|1-w/z|}\\
&=&\frac{\|\nu_0\|_{\infty}}{\pi}\int\limits_{0}^{1}\int\limits_{0}^{2\pi}\frac{rdr\,d\theta}{\rho^2|1-re^{i\theta}/z|}\\
&=&\frac{\|\nu_0\|_{\infty}}{\pi}\int\limits_{0}^{1}\int\limits_{0}^{2\pi}\frac{rdr\,d\theta}{\rho^2|1-re^{i\theta}/\rho|}\\
&=&\frac{\|\nu_0\|_{\infty}}{\pi}\int\limits_{0}^{1}\int\limits_{0}^{2\pi}
\frac{rdr\,d\theta}{\rho^2\sqrt{1+\frac{r^2}{\rho^2}-2\frac{r}{\rho}\cos\,\theta}}\\
&=&\frac{\|\nu_0\|_{\infty}}{\pi}\int\limits_{0}^{1/\rho}\int\limits_{0}^{2\pi}\frac{sds\,d\theta}{\sqrt{1+s^2-2s\cos\,\theta}}\\
&\leq&
\frac{\|\nu_0\|_{\infty}}{\pi}\int\limits_{0}^{1}\int\limits_{0}^{2\pi}\frac{sds\,d\theta}
{\sqrt{1+s^2-2s\cos\,\theta}}\\
&=&\frac{4\|\nu_0\|_{\infty}}{\pi}\int\limits_0^1s\K(s)ds<1.
\end{eqnarray*}
Then $\R p(z,t)>0$ that implies the statement of the theorem.
\end{proof}

 \s

\noindent {\it Remark.}
 If $\|\nu_0 (\cdot,t)\|_{\infty}\leq q$, then
$$\frac{1+|\mu(\zeta,t)|}{1-|\mu(\zeta,t)|}\leq
e^{2tq}\frac{1+|\mu(\zeta,0)|}{1-|\mu(\zeta,0)|}.$$

\s

This obviously follows from the inequality
\[
\frac{\partial |\mu_{f}|}{\partial t}=\frac{\partial
|\mu_{f^{-1}}|}{\partial t}\leq |\dot{\mu}_{f^{-1}}|.
\]

\s

\noindent {\it Remark.} Let us remark that the function $\nu_0$ can
be unilateraly discontinuous on $S^1$ in $\overline{U}$, therefore,
it is not possible, in general, to use the Borel-Pompeiu
formula\index{Borel-Pompeiu formula} to reduce the integral in $p$
to a contour integral.

\s

The equation (\ref{aa7}) is just the L\"owner-Kufarev
equation\index{L\"owner-Kufarev equation}
 in partial derivatives with a special function
$p(z,t)$ given in the above theorems.  \s

Now we discuss the possibility of extending  the equation
(\ref{aa7}) to all of $\mathbb C$. We differentiate the function
$\Phi(\zeta,t,\tau)$ with respect to $\tau$ when $\zeta\in U\cup
U^*$. It follows that
$$
\frac{\partial \Phi(\zeta,t,\tau)}{\partial \tau}\bigg|_{\tau=0}=
\frac{-\overline{f_{\zeta}}}{|f_{\zeta}|^2-|f_{\bar{\zeta}}|^2}\dot{f}+
\frac{{f_{\bar{\zeta}}}}{|f_{\zeta}|^2-|f_{\bar{\zeta}}|^2}\overline{\dot{f}}=:G(\zeta,t).
$$
This formula can be rewritten in the following form
\[
\dot{f}(\zeta,t)=-(f_{\zeta}G(\zeta,t)+f_{\bar{\zeta}}\bar{G}(\zeta,t)).
\]
Taking into account the equation (\ref{aa7}) in $U^*$ we have in the
whole plane
\begin{equation}
\dot{f}(\zeta,t)=\left\{
\begin{array}{ll}
-(f_{\zeta}G(\zeta,t)+f_{\bar{\zeta}}\bar{G}(\zeta,t)),& \mbox{for $\zeta\in \overline{U}$,}\\
-\zeta f_{\zeta} p(\zeta,t), & \mbox{for $\zeta\in U^*$.}
\end{array}\right.\label{aaqq}
\end{equation}
where $p(\zeta,t)$ is a holomorphic in $U^*$ function with the
positive real part by Theorem \ref{quasi}.

The variational formula (\ref{aa5}) and differentiation of the
singular integral  imply that
$G_{\bar{\zeta}}(\zeta,t)=\nu_0(\zeta,t)$, $\zeta\in U$. Now let us
clarify what is $G$. Let us consider $\zeta\in U$. The Pompeiu
formula leads to
\[
G(\zeta,t)=h(\zeta,t)-\frac{1}{\pi}\iint\limits_{U}\frac{\nu_0(w,t)d\sigma_{w}}{w-\zeta},\quad\zeta\in
U
\]
where $h(\zeta,t)$ is a holomorphic function with respect to
$\zeta$. The function $G$ is continuous in $\overline{U}$ and by the
Cauchy theorem
\[
h(\zeta,t)=\frac{1}{2\pi
i}\int\limits_{S^1}\frac{G(w,t)}{w-\zeta}dw.
\]
To obtain the boundary values of the function $G(w,t)$, $|w|=1$, we
will use the second line in (\ref{aaqq}). Unfortunately, in general,
it is not possible to use the same function $f$ in both lines of
(\ref{aaqq}) to obtain boundary values of $G$. Indeed, the mapping
$f(\zeta,t)$ is differentiable regarding to $t$ a.e. in $t\in
[0,t_0)$ locally uniformly in $\zeta \in \mathbb C$, and continuous
in $\zeta\in\mathbb C$ for almost all $t\in [0,t_0)$. Therefore, the
function $-(f_{\zeta}G(\zeta,t)+f_{\bar{\zeta}}\bar{G}(\zeta,t))$,
$\zeta\in U$ is the extension of $-\zeta f'(\zeta,t) p(\zeta,t)$,
$\zeta\in U^*$, whereas $f_{\zeta}$, $\zeta\in U$ is not necessarily
an extension of $f'$, $\zeta\in U^*$.

A simple example of this situation is as follows. Let us consider
the function
\[
f(\zeta,t)=\left\{
\begin{array}{ll}
e^{-t}\left(c\zeta+\dfrac{ \bar{\zeta}}{c}\right),& \mbox{for $\zeta\in \overline{U}$,}\\
e^{-t}\left(c\zeta+\dfrac{1}{c\zeta}\right), & \mbox{for $\zeta\in
U^*$,}
\end{array}\right.
\]
where $c>1$. This mapping forms a subordination chain with the
dilatation $\mu(\zeta)$ that vanishes in $U^*$ and is the constant
$1/c^2$ in $U$. This chain is trivial, but it is not important for
our particular goal here because we do not use at this stage the
crucial variation.  Then,
\[
G(\zeta,t)=\left\{
\begin{array}{ll}
\zeta,& \mbox{for $\zeta\in U$,}\\
\zeta\dfrac{c^2\zeta^2+1}{c^2\zeta^2-1}, & \mbox{for $\zeta\in
U^*$,}
\end{array}\right.
\]
and it splits into two parts that can not be glued on $S^1$. The
same is for the derivatives $f_{\zeta}$ in $U$ and $f'$ in $U^*$.

If $\mu(\zeta,t)$ satisfies the condition (v) in a neighbourhood of
$S^1$ in $\overline{U}$, then the derivatives $f_{\zeta}$,
$f_{\bar{\zeta}}$, $\zeta\in U$ has a continuation onto $S^1$ and
\[
F(\zeta,t)=\dfrac{\overline{f_{\zeta}}\zeta
f'p(\zeta,t)-f_{\bar{\zeta}}\overline{\zeta
f'p(\zeta,t)}}{|f_{\zeta}|^2-|f_{\bar{\zeta}}|^2},\quad \zeta\in
S^1,
\]
where $\zeta f'p(\zeta,t)$ is thought of as the angular limits that
exist a.e. on $S^1$. Moreover, in a neighbourhood of $S^1$ the
derivative $f_{\bar{\zeta}}$ vanishes and
 the function $F(\zeta,t)$ can be written on $S^1$ as
$F(\zeta,t)=\zeta p(\zeta,t)$. In turn,
\[
h(\zeta,t)=\frac{1}{2\pi
i}\int\limits_{S^1}\frac{wp(w,t)}{w-\zeta}dw.
\]

This information allows us formulate the following theorem.

\begin{theorem} Let $f(\zeta,t)$ be a subordination non-trivial chain of
maps in $U^*$ that exists for $t\in [0,t_0)$  and satisfies the
conditions (i--v).

\begin{itemize}
\item[(i)] \ For $\zeta\in U^*$ there exists a holomorphic
function $p(\zeta,t)$ given by Theorem \ref{quasi} such that
\[
\dot{f}(\zeta,t)=-\zeta f'(\zeta,t)p(\zeta,t).
\]

\item[(ii)] \ \ For $\zeta\in U$ there exists a continuous in
$\zeta$ function $F(\zeta,t)$ given by
\[
F(\zeta,t)=\frac{1}{2\pi
i}\int\limits_{S^1}\frac{wp(w,t)}{w-\zeta}dw-
\frac{1}{\pi}\iint\limits_{U}\frac{\nu_0(w,t)d\sigma_{w}}{w-\zeta},
\]
\[
\nu_0(\zeta,t)=\frac{f^{-1}_w}{\bar{f}^{-1}_{\bar{w}}}\frac{\frac{\displaystyle\partial
\mu_{f^{-1}}}{\displaystyle\partial
t}}{1-|\mu_{f^{-1}}|^2}\circ{f(\zeta,t)},
\]
such that
\[
\dot{f}(\zeta,t)=-f_{\zeta}F(\zeta,t)-f_{\bar{\zeta}}\bar{F}(\zeta,t).
\]
\end{itemize}
\end{theorem}

\subsection{The  L\"owner-Kufarev ordinary differential equation}

Dually to the L\"owner-Kufarev partial derivative equation there is
the L\"owner-Kufarev ordinary differential
equation\index{L\"owner-Kufarev equation}. A function $g\in
\Sigma_0$ is represented as a limit
\begin{equation}
\lim\limits_{t\to\infty}e^{-t}w(\zeta,t),\label{lim}
\end{equation}
where the function $w=g(\zeta,t)$ is a solution of the equation
\begin{equation}
\frac{dw}{dt}=-wp(w,t),\label{LKord}
\end{equation}
almost everywhere in $t\in [0,\infty)$, with the initial condition
$g(\zeta,0)=\zeta$. The function $p(\zeta,t)=1+p_1(t)/\zeta+\dots$
is analytic in $U^*$, measurable with respect to $t\in [0,\infty)$,
and its real part $\R\,p(\zeta,t)$ is positive for almost all $t\in
[0,\infty)$. The equation (\ref{LKord}) is known as the {\it
L\"owner-Kufarev ordinary differential equation}. The solutions to
(\ref{LKord}) form a retracting subordination chain $g(\zeta,t)$,
i.e., it satisfies the condition $g(U^*,t)\subset U^*$,
$g(U^*,t)\subset g(U^*,s)$ for $t>s$, and $g(\zeta,0)\equiv \zeta$.

 The connection between (\ref{aa7}) and
(\ref{LKord}) can be thought of as follows. Solving (\ref{aa7}) by
the method of characteristics and assuming $s$ as the parameter
along the characteristics we have $$ \frac{dt}{ds}=1,\quad
\frac{d\zeta}{ds}=\zeta p(\zeta,t), \quad \frac{df}{ds}=0,$$ with
the initial conditions $t(0)=0$, $\zeta(0)=\zeta_0$,
$f(\zeta,0)=f_0(\zeta)$, where $\zeta_0$ is in $U^*$. We see that
the equation (\ref{LKord}) is exactly the characteristic equation
for (\ref{aa7}). Unfortunately, this approach requires the extension
of $f_0(w^{-1}(\zeta,t))$ into $U^*$ because the solution of the
function $f(\zeta,t)$ is given as $f_0(w^{-1}(\zeta,t))$, where
$\zeta=w(\zeta_0,s)$ is the solution of the initial value problem
for the characteristic equation.

Our goal is to deduce a form of the function $p$ on the case of the
subclass $\Sigma_0^{qc}$. Let a one-parameter family of maps
$w=g(\zeta,t)$, $g\in \Sigma_0^{qc}$, satisfy the following
conditions.

\begin{itemize}

\item[(i)] \ \ The maps $g(\zeta,t)$ form a retracting
subordination chain $g(U^*,0)\subset U^*$.

\item[(ii)] \ \ \  The map $g(\zeta,t)$ is meromorphic in $U^*$,
$f(\zeta,t)=\alpha(t)\zeta+a_0(t)+a_1(t)/\zeta+\dots$, where
$\alpha(t)>0$ and differentiable with respect to $t$.

\item[(iii)] \ \ \ The map $g(\zeta,t)$ is a quasiconformal
homeomorphism of $\overline{\mathbb C}$.

\item[(iv)] \ \ \ The chain of maps $g(\zeta,t)$ is not trivial.

\item[(v)] \ \ \ The Beltrami coefficient $\mu_g(\zeta,t)$ of this
map is differentiable with respect to $t$ locally uniformly in $U$.
\end{itemize}

Note that in this case we need not a strong assumption (v) in
Section 6.5.2.

 Set
$$H(\zeta,t,\tau)=g(g(\zeta,t),\tau)=\beta(\tau)w+b_0(\tau)+\frac{b_1(\tau)}{w}+\dots,$$ where $w=g(\zeta,t)$. For
each fixed $t$ the mapping $g(\zeta,t)$ generates a smooth semi-flow
$H^{\tau}$ in $\mathcal G^{qc}$ which is tangent to the path
$g(\zeta,t+\tau)$ at $\tau=0$. Therefore, we use the velocity vector
$wp(w,t)$ (that may be only measurable regarding to $t$) with
$w=g(\zeta,t)$ and obtain $$\frac{\partial
H(\zeta,t,\tau)}{\partial\tau}\Big|_{\tau=0}=g(\zeta,t)p(g(\zeta,t),t).$$
As before, the trajectory $g(\zeta,t)$ generates a pencil of tangent
smooth semi-flows with the tangent vectors $wp(w,t)$, $w=g(z,t)$.
 Since $g(U^*,t)\in
U^*$ for any $t>0$, we can consider the limit
$$\lim\limits_{\tau\to 0}\frac{H(\zeta, t,\tau)-g(\zeta,t)}{\tau
g(\zeta,t)}.$$ We have that
\begin{equation}
\frac{\partial H(\zeta,t,\tau)}{\partial
\tau}\Big|_{\tau=0}=\frac{\partial g(\zeta,t)}{\partial
t}=g(\zeta,t)p(g(\zeta,t),t),\label{aa9}
\end{equation}
 where
$p(\zeta,t)=p_0(t)+p_1(t)/\zeta+\dots$ is an analytic function in
$U^*$ that has  positive real part for almost all fixed $t$. The
equation defined by (\ref{aa9}) is an evolution equation for the
path $g(\zeta,t)$ and the initial condition is given by
$g(\zeta,0)=\zeta$.

We suppose that all $g(\zeta,t)$ admit real analytic quasiconformal
extensions and the family is non-trivial in the above sense. The
function $g(w,\tau)=(H(\zeta,t,\tau)-b_0(\tau))/\beta(\tau)$ can be
extended to a function from $\Sigma_0^{qc}$ and it represents an
equivalence class $[g^{\tau}]\in T$. There is a one-parameter path
$y^{\tau}\in T$ that corresponds to a tangent velocity vector
$\nu(w,t)$ such that
$$\mu_g(w,\tau)=\tau\nu(w,t)+o(\tau),\quad w=g(z,t). $$ We
calculate explicitly the velocity vector making use of the Beltrami
coefficient for a superposition:
\[
\nu(w,t)=\lim\limits_{\tau\to 0}\frac{\mu_{g(w,\tau)}\circ
g(\zeta,t)}{\tau}=\lim\limits_{\tau\to
0}\frac{1}{\tau}\frac{\mu_{H(\zeta,t,\tau)}-\mu_{g(\zeta,t)}}{1-\bar{\mu}_{g(\zeta,t)}\mu_{H(\zeta,t,\tau)}}
\frac{g_{\zeta}(\zeta,t)}{\bar{g}_{\bar{\zeta}}(\zeta,t)},
\]
or
\begin{equation}
\nu(w,t)=\frac{\frac{\partial \mu_{g(\zeta,t)}}{\partial
t}}{1-|\mu_{g(\zeta,t)}|^2}\frac{g_{\zeta}}{\bar{g}_{\bar{\zeta}}}\circ
g^{-1}(w,t), \quad \zeta\in U.\label{aa10}
\end{equation}
It is natural to implement an intrinsic parametrization using the
Teichm\"uller distance $\tau_T(0,[g^t])=t$, and assume the conformal
radius to be $\beta(t)=e^{t}$ that implies $p_0=1$. The assumption
of non-triviality allows us to use the variational formula
(\ref{aa5}) to state the following theorem.

\begin{theorem} Let $g(\zeta,t)$ be a retracting subordination chain of maps
defined in $t\in [0,t_0)$ and $\zeta\in U^*$. Each $g(\zeta,t)$ is a
homeomorphism of $\overline{\mathbb C}$ which is meromorphic in
$U^*$, $\displaystyle g(\zeta,t)=e^{t}\zeta+b_1/\zeta+\dots$, with a
$e^{2t}$-quasiconformal extension to $U$ given by a Beltrami
coefficient $\mu(\zeta,t)$ that is differentiable regarding to $t$
a.e. in $[0, t_0)$. The initial condition is $g(\zeta,0)\equiv
\zeta$.
 Then,  there is a function $p(\zeta,t)$ such that that $\R p(\zeta,t)>0$ for $\zeta\in
 U^*$, and
$$
p(w,t)=1-\frac{1}{\pi}\iint\limits_{g(U,t)}\frac{\nu(u,t)d\sigma_{u}}{w(u-w)},
\quad w\in g(U^*,t),$$ where $\nu(u,t)$ is given by the formula
(\ref{aa10}), $\|\nu\|_{\infty}<1$,
 and $w=g(\zeta,t)$ is a solution to the differential equation
\begin{equation}
\frac{d w}{d t}=w p(w,t),\quad w\in g(U^*,t),\label{aaLK}
\end{equation}
 with the
initial condition $g(\zeta,0)=\zeta$.
\end{theorem}

\s

\noindent {\it Remark.}  Taking into account the superposition we
have
\[
p(g(\zeta,t),t)=1-\frac{1}{\pi}\iint\limits_{U}\frac{
\dot{\mu}_{g}g_{u}^2(u,t)d\sigma_{u}}{g(\zeta,t)(g(u,t)-g(\zeta,t))},
\]
where $u\in U$, $\zeta\in U^*$.

\s

\noindent {\it Remark.} The function $wp(w,t)$ has a continuation
into $g(U,t)$ given by
\[
\frac{d w}{d t}=F(w,t),
\]
where the function $F(w,t)$ is a solution to the equation
\[
\frac{\partial F}{\partial
\bar{w}}=\frac{g_{\zeta}^2\dot{\mu}_g}{|g_{\zeta}|^2-|g_{\bar{\zeta}}|^2}\circ
g^{-1}(w,t).
\]
In contrary to the L\"owner-Kufarev equation in partial derivatives,
the function $F$ is the continuation of $p$ in $U$ through $S^1$.
The solution exists by the Pompeiu integral and can be written as
\begin{eqnarray*}
F(w,t)&=&h(w,t)-\frac{1}{\pi}\iint\limits_{g(U,t)}\frac{g_{\zeta}^2\dot{\mu}_g}{|g_{\zeta}|^2-|g_{\bar{\zeta}}|^2}\circ
g^{-1}(u,t)\frac{d\,\sigma_{u}}{u-w}\\
&=&h(w,t)-\frac{1}{\pi}\iint\limits_{g(U,t)}\frac{\nu(u,t)d\,\sigma_{u}}{u-w},
\end{eqnarray*}
where $w\in g(U,t)$, $h(w,t)$ is a holomorphic functions with
respect to $w$, that can be written as
\[
h(w,t)=\frac{1}{2\pi i}\int\limits_{\partial g(U,t)}\frac{u
p(u,t)}{u-w}du.
\]

Reciprocally, given a function $F(u,t)$, $u\in g(U,t)$, we can write
the function $p(w,t)$ as
\[
p(w,t)=1-\frac{1}{\pi}\iint\limits_{g(U,t)}\frac{F_{\bar{u}}(u,t)d\sigma_{u}}{w(u-w)}
,
\]
where $w\in g(U^*,t)$.

\subsection{Univalent functions smooth on the boundary}

Let us consider the class $\tilde{\Sigma}$ of functions
$f(\zeta)=\alpha\zeta+a_0+a_1/\zeta+\dots$, $\zeta\in U^*$, such
that being extended onto $S^1$ they are $C^{\infty}$ on $S^1$.
Repeating considerations of the preceding subsection for the
embedding of $M$ into the Teichm\"uller space $T$ we come to the
following theorem.

\begin{theorem}\label{smooth} Let $f(\zeta,t)$ be a non-trivial subordination chain of maps
that exists for $t\in [0,t_0)$ and $\zeta\in U^*$. Each $f(\zeta,t)$
is a homeomorphism  $U^*\to \overline{\mathbb C}$ and belongs to
$\tilde{\Sigma}$ for every fixed $t$. All these maps have
quasiconformal extensions to $U$ and there are a real-valued
function $p_0(t)>0$, complex-valued functions $p_1(t)$, real-valued
$C^{\infty}$ functions $d(e^{i\theta},t)$  such that $\R
p(\zeta,t)>0$ for $\zeta\in U^*$,
$$p(\zeta,t)=p_0(t)+\frac{p_1(t)}{\zeta}-\frac{1}{2\pi}\int\limits_{0}^{2\pi}
\frac{e^{i2\theta}d(e^{i\theta},t)d\theta}{\zeta(e^{i\theta}-\zeta)},
\quad \zeta\in U^*,$$  and $f(\zeta,t)$ satisfies the differential
equation $$\frac{\partial f(\zeta,t)}{\partial t}=-\zeta
\frac{\partial f(\zeta,t)}{\partial \zeta}p(\zeta,t),\quad \zeta\in
U^*. $$
\end{theorem}

Theorems \ref{quasi} and \ref{smooth} are linked as follows. For a
given subordination chain of maps $f(\zeta,t)\in \tilde{\Sigma}$,
that exists for $t\in [0,t_0)$ and $\zeta\in U^*$, there is a
$C^{\infty}$ function $d(e^{i\theta},t)$ by Theorem \ref{smooth} and
we can construct the function $\nu(\zeta,t)$ by the Douady-Earle
extension and the formula (\ref{mu1}). Then, the function
 $f(\zeta,t)$ satisfies the equation of Theorem \ref{quasi} with
 $p(\zeta,t)$ defined by such $\nu(\zeta,t)$.

\s

Let us consider the ordinary L\"owner-Kufarev
equation\index{L\"owner-Kufarev equation} for the functions smooth
on $S^1$. If the retracting chain $g(\zeta,t)$ is smooth on $S^1$,
then we use again the embedding of $M$ into $T$ and reach a similar
result.

\begin{theorem} Let $g(\zeta,t)$ be a retracting non-trivial subordination chain of normalized maps
that exists for $t\in [0,t_0)$ and $\zeta\in U^*$. Each $g(\zeta,t)$
is meromorphic in $U^*$, smooth on $S^1$, and $\displaystyle
g(\zeta,t)=\beta(t)\zeta+b_0(t)+\frac{b_1(t)}{\zeta}+\dots$,
$\beta(t)>0$. An additional assumption is that $g:\, U^*\to U^*$ for
each fixed $t$.
 Then, there are a real-valued function $p_0(t)$, a complex-valued function
 $p_1(t)$, and
  a smooth real-valued function $d(e^{i\theta},t)$, such that
$\R p(\zeta,t)>0$ for $\zeta\in U^*$,
$$p(\zeta,t)=p_0(t)+\frac{p_1(t)}{\zeta}-\frac{1}{2\pi
i}\int\limits_{S^1}\left(\frac{zg'(z,t)}{g(z,t)}\right)^2\frac{d
(z,t)dz}{g(z,t)-\zeta}, \quad \zeta\in U^*,$$
 and $w=g(\zeta,t)$ is a solution to the differential equation
$$\frac{d w}{d t}=w p(w,t),\quad w\in g(U^*,t) $$ with the
initial condition $g(\zeta,0)=\zeta$.
\end{theorem}

\s

\noindent {\it Remark.} If we work with normalized functions
$$g(\zeta,t)=e^t\zeta+\frac{b_1(t)}{\zeta}+\dots,$$
then $p_0(t)\equiv 1$, $p_1(t)\equiv 0$.

\subsection{An application to Hele-Shaw flows}

Theorem \ref{smooth} is linked to the Hele-Shaw\index{Hele-Shaw
problem} free boundary problem as follows. Starting with a smooth
boundary $\Gamma_0$ the one-parameter family $\Gamma(t)$ consists of
smooth curves as long as the solutions exist. Let us consider the
equation (\ref{PGL}). Under injection we have a subordination chain
of domains $\Omega(t)$. The Schwarz kernel can be developed as
$$\frac{\zeta+e^{i\theta}}{\zeta-e^{i\theta}}=1+\frac{2e^{i\theta}}{\zeta}+
\frac{2e^{2i\theta}}{\zeta(\zeta-e^{i\theta})}. $$ Therefore, in
Theorem \ref {smooth} we can put
$$p_0(t)=\frac{1}{2\pi}\int\limits_{0}^{2\pi}\frac{1}{|f'(e^{i\theta},t)|^2}d\theta,
\quad
p_1(t)=\frac{1}{\pi}\int\limits_{0}^{2\pi}\frac{e^{i\theta}}{|f'(e^{i\theta},t)|^2}d\theta,$$
and $$d(e^{i\theta},t)=\frac{-2}{|f'(e^{i\theta},t)|^2}.$$

Apart from the trivial elliptic case there are no self-similar
solutions, and therefore the Hele-Shaw dynamics $f(\zeta,t)$
generates a non-trivial path in $T$. Thus, given a Hele-Shaw
evolution $\Gamma(t)=f(S^1,t)$ we observe a differentiable
non-trivial path on $T$, such that at any time $t$ the tangent
vector $\nu$ is a harmonic Beltrami differential given by
\[
\nu(\zeta,t)=\frac{-3}{\pi}\int\limits_0^{2\pi}\frac{(1-|\zeta|^2)^2}{(1-e^{i\theta}\bar{\zeta})^4}
\frac{e^{2i\theta}}{|f'(e^{i\theta},t)|^2}d\theta.
\]
The corresponding co-tangent vector is
\[
\varphi(\zeta,t)=\frac{6}{\pi}\int\limits_0^{2\pi}
\frac{e^{-2i\theta}d\theta}{(1-e^{-i\theta}\zeta)^4|f'(e^{i\theta},t)|^2}.
\]


\begin{thebibliography}{99}

\bibitem{Alexandrov}
I.~A.~Aleksandrov, {\it Parametric continuations in the theory of
univalent functions}, Nauka, Moscow, 1976. (in Russian)

\bibitem{Becker}
J.~Becker, {\it L\"ownersche Differentialgleichung und
quasikonform fortsetzbare schlichte Funktionen}, J. Reine Angew.
Math. {\bf 255} (1972), 23--43.

\bibitem{Becker1}
J.~Becker, {\it L\"ownersche Differentialgleichung und
Schlichtheitskriterien}, Math. Ann. {\bf 202} (1973), 321--335.

\bibitem{Becker2}
J.~Becker, {\it Conformal mappings with quasiconformal
extensions}. Aspects of contemporary complex analysis (Proc. NATO
Adv. Study Inst., Univ. Durham, Durham, 1979) Academic Press,
London-New York, 1980, 37--77.

\bibitem{BeckPom}
J.~Becker, Ch.~Pommerenke, {\it On the Hausdorff dimension of
quasicircles}, Ann. Acad. Sci. Fenn. Ser. A I Math. {\bf 12}
(1987), no. 2, 329--333.


\bibitem{BlAhl} A.~Beurling, L.~V.~Ahlfors, {\it
The boundary correspondence under quasiconformal mappings}, Acta
Math. {\bf 96} (1956), 125--142.

\bibitem{Bib}
L.~Bieberbach, {\it \"Uber die Koeffizienten derjenigen
Potenzreihen, welche eine schlichte Abbildung des Einheitskreises
vermitteln}, S.-B. Preuss. Akad. Wiss. (1916), 940--955.


\bibitem{Branges}
L.~de~Branges, {\it A proof of the Bieberbach conjecture}, Acta
Math. {\bf 154} (1985) no. 1--2, 137--152.


\bibitem{DE}
A.~Douady, C.~J.~Earle, {\it Conformally natural extension of
homeomorphisms of the circle}, Acta Math. {\bf 157} (1986), no.
1-2, 23--48.

\bibitem{Duren}
P.~Duren, {\it Univalent functions}, Springer, New York, 1983.

\bibitem{Galin} L.~A.~Galin, {\it Unsteady filtration with a free
surface}, Dokl. Akad. Nauk USSR, {\bf 47} (1945), 246--249. (in
Russian)

\bibitem{GL}
F.~Gardiner, N.~Lakic, {\it Quasiconformal Teichm\"uller theory},
Mathematicsl Surveys and Monographs, vol. 76, Amer. Math. Soc.,
2000.

\bibitem{Gehring}
F.~W.~Gehring, E.~Reich, {\it Area distortion under quasiconformal
mappings}, Ann. Acad. Sci. Fenn. Ser. A I, No. 388, 1966, 15 pp.

\bibitem{Gill} K.~A.~Gillow, S.~D.~Howison, {\it A bibliography of free and
moving boundary problems for Hele-Shaw and Stokes flow}, Published
electronically at
URL$>$http://www.maths.ox.ac.uk/\~{}howison/Hele-Shaw.

\bibitem{Gor}
V.~V.~Goryainov, {\it Fractional iterates of functions that are
analytic in the unit disk with given fixed points}, Mat. Sb. {\bf
182} (1991), no. 9, 1281--1299; Engl. Transl. in Math. USSR-Sb.
{\bf 74} (1993), no. 1, 29--46.

\bibitem{Gor1}
V.~V.~Goryainov, {\it One-parameter semigroups of analytic
functions}, Geometric function theory and applications of complex
analysis to mechanics: studies in complex analysis and its
applications to partial differential equations, 2 (Halle, 1988),
Pitman Res. Notes Math. Ser., 257, Longman Sci. Tech., Harlow,
1991, 160--164.

\bibitem{Gor2}
V.~V.~Goryainov, {\it One-parameter semigroups of analytic
functions and a compositional analogue of infinite divisibility.}
Proceedings of the Institute of Applied Mathematics and Mechanics,
Vol. 5,  Tr. Inst. Prikl. Mat. Mekh., 5,
   Nats. Akad. Nauk Ukrainy Inst. Prikl. Mat. Mekh., Donetsk, 2000,
   44--57.(in Russian)

\bibitem{Bjorn} B.~Gustafsson, {\it On a differential equation arising in a Hele-Shaw
flow moving boundary problem}, Arkiv f\"or Mat., {\bf 22} (1984),
no. 1,  251--268.

\bibitem{Gut}
V.~Ya.~Gutlyanski\u\i, {\it Parametric representation of univalent
functions}, Dokl. Akad. Nauk SSSR {\bf 194} (1970), 750--753;
Engl. Transl. in Soviet Math. Dokl. {\bf 11} (1970), 1273--1276.

\bibitem{Gut2}
V.~Ya.~Gutlyanski\u\i, {\it On some classes of univalent
functions}.  Theory of functions and mappings, "Naukova Dumka",
Kiev, 1979,  85--97. (in Russian)

\bibitem{Gut3}
V.~Ya.~Gutljanski\u\i, {\it The method of variations for univalent
analytic functions with a quasiconformal extension}, Sibirsk. Mat.
Zh. {\bf 21} (1980), no. 2, 61--78; Engl. Transl. in Siberian
Math. J. {\bf 21} (1980), no. 2, 190--204.

\bibitem{He} Cheng Qi He, {\it A parametric representation of quasiconformal
extensions}, Kexue Tongbao, {\bf 25} (1980), no. 9, 721--724.


\bibitem{Hele} H.~S.~Hele-Shaw, {\it The flow of water},
Nature {\bf 58} (1898), no. 1489, 33--36.

\bibitem{How}
S.~D.~Howison, {\it Complex variable methods in Hele-Shaw moving
boundary problems},  European J. Appl. Math., {\bf 3} (1992), no.
3, 209--224.

\bibitem{KY1}
A.~A.~Kirillov, D.~V.~Yuriev, {\it K\"ahler geometry of the
infinite-dimensional homogeneous space $M={\rm Diff}\sb +(S\sp
1)/{\rm Rot}(S\sp 1)$}, Funktsional. Anal. i Prilozhen. {\bf 21}
(1987), no. 4, 35--46. (in Russian)

\bibitem{KY2}
A.~A.~Kirillov, D.~V.~Yuriev, {\it Representations of the Virasoro
algebra by the orbit method}, J. Geom. Phys. {\bf 5} (1988), no.
3, 351--363.

\bibitem{Kir}
A.~A.~Kirillov, {\it K\"ahler structure on the $K$-orbits of a
group of diffeomorphisms of the circle}, Funktsional. Anal. i
Prilozhen. {\bf 21} (1987), no. 2, 42--45.

\bibitem{Kufarev}
P.~P.~Kufarev, {\it On one-parameter families of analytic
functions}, Rec. Math. [Mat. Sbornik] N.S. {\bf 13(55)} (1943),
87--118.

\bibitem{Lawryn}
J.~{\L}awrynowicz, J.~Krzy{\.{z}}, {\it Quasiconformal mappings in
the plane: parametric methods}, Lecture Notes in Math., Vol. 978,
Springer-Verlag, Berlin-New York, 1983.

\bibitem{Lecko}
A.~Lecko, D.~Partyka, {\it An alternative proof of a result due to
Douady and Earle}, Annales Univ. Mariae Curie-Sk{\l}odowska, Sect.
A, {\bf 42}  (1988), 59--68.

\bibitem{Lehto}
O.~Lehto, {\it Univalent functions and Teichm\"uller spaces},
Graduate Texts in Mathematics, 109. Springer-Verlag, New York, 1987.

\bibitem{Loewner}
K.~L\"owner, {\it Untersuchungen \"uber schlichte konforme
Abbildungen des Einheitskreises}, Math. Ann. {\bf 89} (1923),
103--121.

\bibitem{Nagbook} S.~Nag, {\it The complex analytic theory of Teichm\"uller
spaces}, Wiley-Interscience Publ. John Wiley \& Sons, Inc., New
York, 1988.

\bibitem{Nehari}
Z.~Nehari, {\it Schwarzian derivatives and schlicht functions},
Bull. Amer. Math. Soc. {\bf 55} (1949), no. 6, 545--551.

\bibitem{Ock}
H.~Ockendon, J.~R.~Ockendon, {\it Viscous Flow}, Cambridge U.P.,
1995.

\bibitem{Polub1} P.~Ya.~Polubarinova-Kochina, {\it On a problem of the
motion of the contour of a petroleum shell}, Dokl. Akad. Nauk
USSR, {\bf 47} (1945), no. 4, 254--257. (in Russian)

\bibitem{Polub2} P.~Ya.~Polubarinova-Kochina, {\it Concerning unsteady
motions in the theory of filtration}, Prikl. Matem. Mech., {\bf 9}
(1945), no. 1, 79--90. (in Russian)

\bibitem{Pom1}
Ch.~Pommerenke, {\it \"Uber die Subordination analytischer
Funktionen}, J. Reine Angew. Math. {\bf 218} (1965), 159--173.

\bibitem{Pom2}
Ch.~Pommerenke, {\it Univalent functions, with a chapter on
quadratic differentials by G.~Jensen}, Vandenhoeck \& Ruprecht,
G\"ottingen, 1975.

\bibitem{Reissig}
M.~Reissig, L.~Von Wolfersdorf, {\it A simplified proof for a
moving boundary problem for Hele-Shaw flows in the plane},  Ark.
Mat., {\bf 31} (1993), no. 1,  101--116.

\bibitem{Saffman} P.~G.~Saffman, G.~I.~Taylor, {\it The penetration of a fluid into
a porous medium or Hele-Shaw cell containing a more viscous
liquid}, Proc. Royal Soc. London, Ser. A {\bf 245} (1958), no.
281, 312--329.

\bibitem{Shah}
Shah Dao-Shing, {\it Parametric representation of quasiconformal
mappings}, Science Record {\bf 3} (1959), 400--407.

\bibitem{Shoikhet}
D.~Shoikhet, {\it Semigroups in geometrical function theory},
Kluwer Academic Publishers, Dordrecht, 2001.

\bibitem{Petersson} {S.~Wolpert}, {\it Thurston's Riemannian metric for Teichm\"uller
space},  J. Differential Geom., {\bf 23} (1986), no. 2, 143--174.

\bibitem{Vas} A.~Vasil'ev, {\it Univalent functions in two-dimensional free boundary
problems}, Acta Applic. Math. {\bf 79} (2003), no. 3, 249--280

\bibitem{VinKuf}
Yu.~P.~Vinogradov, P.~P.~Kufarev, {\it On a problem of
filtration}, Akad. Nauk SSSR. Prikl. Mat. Meh., {\bf 12} (1948),
181--198. (in Russian)
\end{thebibliography}
\end{document}